\documentclass{amsart}

\usepackage{amsthm}
\usepackage{amsfonts}
\usepackage{amsmath}
\usepackage{amssymb}
\usepackage{mathrsfs}
\usepackage{fullpage}
\usepackage{stmaryrd}
\usepackage{hyperref}
\usepackage{verbatim}

\theoremstyle{plain}

\theoremstyle{definition}

\theoremstyle{remark}

\newcommand{\Begin}[2]{\begin{#1}\label{#2}}

\newcommand{\R}{\mathbb{R}}

\newcommand{\bPi}{\mathbf{\Pi}}
\newcommand{\bSigma}{\mathbf{\Sigma}}
\newcommand{\bDelta}{\mathbf{\Delta}}
\newcommand{\bbC}{\mathbb{C}}

\newcommand{\bbP}{\mathbb{P}}

\newcommand{\bbR}{\mathbb{R}}

\newcommand{\SCRL}{\mathscr{L}}

\newcommand{\forces}{\Vdash}
\newcommand{\analytic}{{\bSigma_1^1}}

\newcommand{\coanalytic}{{\bPi_1^1}}

\newcommand{\borel}{{\bDelta_1^1}}

\newcommand{\bairespace}{{{}^\omega\omega}}

\newcommand{\finNaturalSequence}{{{}^{<\omega}\omega}}

\begin{document}

\title{When an Equivalence Relation with All Borel Classes will be Borel Somewhere?}

\author{William Chan}
\address{Department of Mathematics, California Institute of Technology, Pasadena, CA 91125}
\email{wcchan@caltech.edu}

\author{Menachem Magidor}
\address{Einstein Institute of Mathematics, Hebrew University of Jerusalem}
\email{mensara@savion.huji.ac.il}

\begin{abstract}
In $\mathsf{ZFC}$, if there is a measurable cardinal with infinitely many Woodin cardinals below it, then for every equivalence relation $E \in L(\bbR)$ on $\bbR$ with all $\borel$ classes and every $\sigma$-ideal $I$ on $\bbR$ so that the associated forcing $\bbP_I$ of $I^+$ $\borel$ subsets is proper, there exists some $I^+$ $\borel$ set $C$ so that $E \upharpoonright C$ is a $\borel$ equivalence relation. In $\mathsf{ZF} + \mathsf{DC} + \mathsf{AD}_\bbR + V = L(\mathscr{P}(\bbR))$, for every equivalence relation $E$ on $\bbR$ with all $\borel$ classes and every $\sigma$-ideal $I$ on $\bbR$ so that the associated forcing $\bbP_I$ is proper, there is some $I^+$ $\borel$ set $C$ so that $E \upharpoonright C$ is a $\borel$ equivalence relation.
\end{abstract}

\maketitle\let\thefootnote\relax\footnote{August 17, 2016
\\*\indent The first author was partially supported by NSF grants DMS-1464475 and EMSW21-RTG DMS-1044448}


\section{Introduction}\label{Introduction}
The basic question of interest is: 

\Begin{question}{basic question}
If $E$ is an equivalence relation on $\bairespace$, is $E$ a simpler equivalence relation when restricted to some subset?
\end{question}

This question can also be asked for equivalence relations on arbitrary Polish spaces, but for simplicity, this paper will only consider equivalence relations on $\bairespace$. Usually, descriptive set theoretic results about $\bairespace$ have proofs that can be transfered to arbitrary Polish spaces. 

What should be the measure of complexity and what should be the paragon of simplicity? The measure of complexity will vaguely be definability and there is no need to explicitly state what it is since the paper will only strive to reach the base of complexity. However, there are various useful notions of definability given by considerations in topology, recursion theory, logical complexity, and set theory. The base of definable complexity needs to be explicitly stated. The class of Borel sets (denoted $\borel$) is chosen to be this base since it is a simple class characterized by all the notions of definability mentioned above. Moreover, many natural mathematical concerns appear at this level, and $\borel$ objects seem to be well behaved and relatively well understood. 

Now the question can be more precisely formulated: 

\Begin{question}{basic question version 2}
If $E$ is an equivalence relation on $\bairespace$, is there a $\borel$ set $C \subseteq \bairespace$ so that $E \upharpoonright C$ is a $\borel$ equivalence relation?
\end{question}

Here, $E \upharpoonright C = E \cap (C \times C)$. However, there is one obvious triviality. If $C$ is countable, then any equivalence relation restricted to $C$ is $\borel$. Since countable subsets of $\bairespace$ belong to any $\sigma$-ideal on $\bairespace$ which contains all singletons, this egregious triviality disappears if one asks that, in the above question, $C$ be $\borel$ and non-trivial according to a $\sigma$-ideal on $\bairespace$. Subsets of $\bairespace$ that are not in the ideal $I$ are called $I^+$ sets. In this paper, $\sigma$-ideals will always contain all the singletons.

However, it is unclear how to approach this question for arbitrary $\sigma$-ideals. The collection of available techniques is greatly enriched by considering $\sigma$-ideals on $\bairespace$ so that the associated forcing $\bbP_I$ of $\borel$ $I^+$ sets is a proper forcing. Considering such $\sigma$-ideals makes available powerful tools from models of set theory and absoluteness. (In fact, the questions below all have negative answers when considering arbitrary $\sigma$-ideals. See Section \ref{Basics}.) Now a test question can be posed for a slightly more complicated class of equivalence relations than the $\borel$ equivalence relations: Analytic (denoted $\analytic$) sets are continuous images of $\borel$ or even closed sets.

\Begin{question}{basic question analytic}
Let $E$ be a $\analytic$ equivalence relation on $\bairespace$. Let $I$ be a $\sigma$-ideal on $\bairespace$ so that $\bbP_I$ is a proper forcing. Is there an $I^+$ $\borel$ set $C$ so that $E \upharpoonright C$ is a $\borel$ equivalence relation?
\end{question}

Note that questions like the above are very familiar. For example, the ideal of Lebesgue null set and the ideal of meager sets have the property that their associated forcings are proper forcings. It is very common in mathematics to ask questions about properties that hold on positive measure sets (or Lebesgue almost everywhere) or on non-meager (or comeager) sets. 

Unfortunately, Question \ref{basic question analytic} has a negative answer:

\Begin{proposition}{negative answer analytic}
There is a $\analytic$ equivalence relation $E$ and a $\sigma$-ideal $I$ with $\bbP_I$ proper so that for all $\borel$ $I^+$ set $C$, $E \upharpoonright C$ is not $\borel$. 
\end{proposition}

\begin{proof}
See \cite{Canonical-Ramsey-Theory-on-Polish-Spaces}, Example 4.25.
\end{proof}

So a positive answer is not even possible for the simplest class of equivalence relations in the projective hierarchy just above $\borel$. A positive answer to any variation of the basic question will likely only be feasible if the equivalence relations bear at least some resemblance to $\borel$ equivalence relations. \cite{Canonical-Ramsey-Theory-on-Polish-Spaces} then proved that a positive answer does hold for $\analytic$ equivalence relations with all countable classes and equivalence relations $\borel$ reducible to orbit equivalence relations of Polish groups actions. In both these examples, the equivalence relations have all $\borel$ classes. Of course, $\borel$ equivalence relations have all $\borel$ classes. Perhaps those two examples give evidence that a sufficient resemblance for a positive answer is the property of having all $\borel$ classes. \cite{Canonical-Ramsey-Theory-on-Polish-Spaces} asked the following question:

\Begin{question}{main question analytic}
(\cite{Canonical-Ramsey-Theory-on-Polish-Spaces} Question 4.28) Let $E$ be a $\analytic$ equivalence relation on $\bairespace$ with all $\borel$ classes. Let $I$ be a $\sigma$-ideal on $\bairespace$ so that $\bbP_I$ is a proper forcing. Let $B$ be an $I^+$ $\borel$ set. Is there some $C \subseteq B$ which is $I^+$ $\borel$ so that $E \upharpoonright C$ is a $\borel$ equivalence relation?
\end{question}

Under large cardinal assumptions, this question has a positive answer: Here coanalytic set (denoted $\coanalytic$) are complements of $\analytic$ sets.

\Begin{theorem}{positive answer main question analytic}
Suppose for all $X \in H_{(2^{\aleph_0})^+}$, $X^\sharp$ exists. Then for all $\analytic$ and $\coanalytic$ equivalence relations with all $\borel$ classes, any $\sigma$-ideal $I$ on $\bairespace$ with $\bbP_I$ proper, and $B$ an $I^+$ $\borel$ set, there exists some $I^+$ $\borel$ set $C \subseteq B$ so that $E \upharpoonright C$ is $\borel$.
\end{theorem}

\begin{proof}
See \cite{Equivalence-Relations-Which-Are-Borel-Somewhere}. Also see \cite{Borel-Canonization-of-Analytic-Sets} for a similar result proved using a measurable cardinal.
\end{proof}

It should be noted that the proofs of Theorem \ref{positive answer main question analytic} in both \cite{Equivalence-Relations-Which-Are-Borel-Somewhere} and \cite{Borel-Canonization-of-Analytic-Sets} use an approximation of $\analytic$ equivalence relations by $\borel$ equivalence relations: Burgess showed that for every $\analytic$ equivalence relation $E$ there is (in a uniform way) an $\omega_1$-length decreasing sequence $(E_\alpha : \alpha < \omega_1)$ of $\borel$ equivalence relations so that $E = \bigcap_{\alpha < \omega_1} E_\alpha$. The strategy of the proof is to find some countable elementary $M \prec H_\Xi$, where $\Xi$ is large enough to contain certain desired objects, and some countable ordinal $\alpha$ so that if $C$  is the $I^+$ $\borel$ set of $\bbP_I$-generic reals over $M$ (which exists by properness of $\bbP_I$), then $E \upharpoonright C = E_\alpha \upharpoonright C$. The sharps are used to obtain the absoluteness necessary to determine the countable level $\alpha$ at which the $E$ classes and $E_\alpha$ classes of all generic reals stabilize.

In conversation with the first author, Neeman asked the following generalization of Question \ref{main question analytic}: Projective sets are those obtainable by applying finitely many applications of complements and continuous images starting with the $\borel$ sets.

\Begin{question}{main question projective}
Assume some large cardinal hypotheses. Let $E$ be a projective equivalence relation with all $\borel$ classes. Let $I$ be a $\sigma$-ideal on $\bairespace$ with $\bbP_I$ proper. Let $B \subseteq \bairespace$ be an $I^+$ $\borel$ subset. Does there exists some $I^+$ $\borel$ $C \subseteq B$ so that $E \upharpoonright C$ is $\borel$?
\end{question}

It is unclear if the proofs of Theorem \ref{positive answer main question analytic} can be generalized to give an answer to this question since there does not appear to be any form of $\borel$ approximation to arbitrary projective equivalence relations. Moreover, it is known to be consistent that there is a negative answer to Question \ref{main question projective} even when restricted to the next level of the projective hierarchy above $\analytic$ and $\coanalytic$. A $\bSigma_2^1$ set is a continuous image of a $\coanalytic$ set. A $\bPi_2^1$ set is the complement of a $\bSigma_2^1$ set. A $\bDelta_2^1$ set is a set that is both $\bSigma_2^1$ and $\bPi_2^1$:

\Begin{proposition}{negative answer delta12 in L}
In the constructible universe $L$, there is a $\bDelta_2^1$ equivalence relation with all classes countable so that for every $\sigma$-ideal $I$ and every $I^+$ $\borel$ set $B$, $E \upharpoonright B$ is not $\borel$.
\end{proposition}

\begin{proof}
See \cite{Equivalence-Relations-Which-Are-Borel-Somewhere} or \cite{Borel-Canonization-of-Analytic-Sets}.
\end{proof}

In fact, it is not even known what is the status of Question \ref{main question analytic} or its $\coanalytic$ analog in $L$. Perhaps the most interesting open question in this area is whether it is consistent that Question \ref{main question analytic} or its $\coanalytic$ analog has a negative answer. See the conclusion section of \cite{Equivalence-Relations-Which-Are-Borel-Somewhere} for some discussions on this question.

This paper will be concerned with extending a positive answer to these types of questions to larger classes of equivalence relations on $\bairespace$ with all $\borel$ classes. As mentioned above, some new methods will need to be developed to take the role of Burgess's approximation in Theorem \ref{positive answer main question analytic}. A certain game will be used to fulfill this role. 

Question \ref{main question projective} will be answered by an even more general result. Like in Theorem \ref{positive answer main question analytic}, the results of this paper will be proved in an extension of $\mathsf{ZFC}$, the standard axioms of set theory. Here, $\mathsf{ZFC}$ will be augmented by large cardinal axioms. The large cardinal axioms used here are well accepted and have proven to be very useful in descriptive set theory.

The model $L(\bbR)$ is the smallest inner model of $\mathsf{ZF}$ (possibly without the axiom of choice) containing all the reals of the original universe. It contains all the sets which are ``constructible'' (in the sense of G\"odel) from the reals of the original universe. Nearly all objects of ordinary mathematics can be found in $L(\bbR)$. In particular, all projective subsets of $\bairespace$ belong to $L(\bbR)$. A main result of the paper is:
\\*
\\*\noindent \textbf{Theorem \ref{L(R) canonicalization}.} Suppose there is a measurable cardinal with infinitely many Woodin cardinals below it. Let $I$ be a $\sigma$-ideal on $\bairespace$ so that $\bbP_I$ is a proper forcing. Let $E \in L(\bbR)$ be an equivalence relation on $\bairespace$. 

If $E$ has all $\analytic$ ($\coanalytic$ or $\borel$) classes, then for every $I^+$ $\borel$ set $B$, there is an $I^+$ $\borel$ $C \subseteq B$ so that $E \upharpoonright C$ is $\analytic$ ($\coanalytic$ or $\borel$, respectively).
\\*
\\*\indent This gives a positive answer to Question \ref{main question projective}. Moreover, it shows that for a large class of equivalence relations on $\bairespace$ so that all the equivalences classes belong to a particular pointclass of the first level of the projective hierarchy, the equivalence relation somewhere is as simple as its equivalence classes.

Having answered Question \ref{main question projective} positively and even given a positive answer for the larger class of $L(\bbR)$ equivalence relation with all $\borel$ classes, the ultimate natural question is the following:

\Begin{question}{every equivalence relation canonicalization}
Is it consistent relative some large cardinals, that (the axiom of choice fails and) for every equivalence relation $E$ with all $\borel$ classes and every $\sigma$-ideal $I$ on $\bairespace$ such that $\bbP_I$ is a proper forcing, there is an $I^+$ $\borel$ set $C$ so that $E \upharpoonright C$ is a $\borel$ equivalence relation?
\end{question}

As it is often the case for various regularity properties like the perfect set property, Lebesgue measurability, or the property of Baire, the axiom of choice can be used with a diagonalization argument to produce a failure of this property. In fact, using the axiom of choice, there is an equivalence relation with classes of size at most two so that for any $\sigma$-ideal $I$ and any $I^+$ $\borel$ set $C$, $E \upharpoonright C$ is not $\borel$.

For the regularity properties mentioned above, it is consistent that all sets have these properties in a choiceless model of $\mathsf{ZF}$, like the model $L(\bbR)$. For instance, if the axiom of determinacy, $\mathsf{AD}$, holds then all sets are Lebesgue measurable and have the property of Baire. 

Assuming determinacy for certain games on the reals, every equivalence relation with all $\borel$ classes can be canonicalized by $\sigma$-ideals whose associated forcings are proper: 
\\*
\\*\noindent \textbf{Theorem \ref{AD(R) canonicalization for all equiv meager ideal}.} Assume $\mathsf{ZF} + \mathsf{DC} + \mathsf{AD}_\bbR $. Let $E$ be an equivalence relation on $\bairespace$. If $E$ has all $\analytic$ ($\coanalytic$ or $\borel$) classes, then for every nonmeager $\borel$ set $B$, there is a $\borel$ set $C \subseteq B$ which is comeager in $B$ so that $E \upharpoonright C$ is $\analytic$ ($\coanalytic$ or $\borel$, respectively).
\\*
\\*\noindent \textbf{Theorem \ref{AD(R) V = L(P(R)) canonicalization for all equiv all proper ideal}.} Assume $\mathsf{ZF} + \mathsf{DC} + \mathsf{AD}_\bbR + V = L(\mathscr{P}(\bbR))$. Let $I$ be a $\sigma$-ideal on $\bairespace$ so that $\bbP_I$ is proper. Let $E$ be an equivalence relation on $\bairespace$. If $E$ has all $\analytic$ ($\coanalytic$ or $\borel$) classes, then for every $I^+$ $\borel$ set $B$, there is an $I^+$ $\borel$ set $C \subseteq B$ so that $E \upharpoonright C$ is $\analytic$ ($\coanalytic$ or $\borel$, respectively).
\\*
\\*\indent Section \ref{Basics} will review the basics of idealized forcing, the theory of measure, and homogeneous trees. The relevant game concepts will be introduced here.

Section \ref{The Game} will prove that certain types of equivalence relations can be $\borel$, $\analytic$, or $\coanalytic$ equivalence relations on $I^+$ $\borel$ subsets of $\bairespace$ for any $\sigma$-ideal $I$ so that $\bbP_I$ is proper, under three assumptions about absoluteness and tree representations. The main results of this section will be proved using a certain game. This section can be understood with just basic knowledge of set theory and forcing. The results of this section holds more generally for relations with $\analytic$ ($\coanalytic$ or $\borel$) sections. Therefore, all the theorems in this paper have an analogous statement for graphs $G$ so that for all $x \in \bairespace$, the set $G_x = \{y : x \ G \ y\}$ is $\analytic$ ($\coanalytic$ or $\borel$). However, this paper will focus mostly on equivalence relations.

Section \ref{Canonicalization} will mostly assume axiom of choice and will give a general situation in which the three assumptions used in the previous section hold. This section will give a very brief survey of the theory of generic absoluteness and tree representations of subsets of $\bairespace$, especially the Martin-Solovay tree construction. Theorem \ref{L(R) canonicalization} will be presented.

Section \ref{Canonicalization for All Equivalence Relations} will assume a bit more than the axiom of determinacy for the reals and will mention the necessary results about tree representations and generic absoluteness to show that the three assumptions from Section \ref{The Game} holds for every equivalence relation with all $\analytic$, $\coanalytic$, or $\borel$ classes. Finally, Theorem \ref{AD(R) canonicalization for all equiv meager ideal} and Theorem \ref{AD(R) V = L(P(R)) canonicalization for all equiv all proper ideal} will be presented.

The authors would like to thank Alexander Kechris, Itay Neeman, and Zach Norwood for many useful discussions about the contents of this paper.

\section{Basics}\label{Basics}
In this paper, $\sigma$-ideals always contain all the singleton.
\Begin{definition}{idealized forcing}
Let $I$ be a $\sigma$-ideal on $\bairespace$. Let $\bbP_I = (\borel \setminus I, \subseteq, \bairespace)$ be the forcing of $I^+$ $\borel$ subsets of $\bairespace$ ordered by $\leq_{\bbP_I} = \subseteq$ and has largest element $1_{\bbP_I} = \bairespace$. Often $\bbP_I$ is identified with $\borel \setminus I$. 
\end{definition}

\Begin{fact}{name for generic element}
Let $I$ be a $\sigma$-ideal on $\bairespace$. There is a name $\dot x_\emph{gen} \in V^{\bbP_I}$ so that for all $\bbP_I$-generic filters $G$ over $V$ and all $\borel$ sets $B$ coded in $V$, $V[G] \models B \in G \Leftrightarrow \dot x_\emph{gen}[G] \in B$. 
\end{fact}

\begin{proof}
See \cite{Forcing-Idealized}, Proposition 2.1.2.
\end{proof}

\Begin{definition}{generic real}
Let $I$ be a $\sigma$-ideal on $\bairespace$. Let $M \prec H_\Xi$ be a countable elementary substructure for some sufficiently large cardinal $\Xi$. $x \in \bairespace$ is $\bbP_I$-generic over $M$ if and only if the collection $\{B \in \bbP_I \cap M : x \in B\}$ is a $\bbP_I$-generic filter over $M$. 
\end{definition}

The following results makes available some very useful techniques for handling ideals whose associated forcings are proper forcings. For the purpose of this paper, the following may as well be taken as the definition of properness:

\Begin{proposition}{properness equivalence}
Let $I$ be a $\sigma$-ideal on $\bairespace$. The following are eqivalent:

(i) $\bbP_I$ is a proper forcing.

(ii) For any sufficiently large cardinal $\Xi$, every $B \in \bbP_I$, and every countable $M \prec H_\Xi$ with $\bbP_I \in M$ and $B \in M$, the set $C = \{x \in B : x \text{ is } \bbP_I\text{-generic over } M\}$ is an $I^+$ $\borel$ set.
\end{proposition}

\begin{proof}
See \cite{Forcing-Idealized}, Proposition 2.2.2.
\end{proof}

This proposition shows that $\sigma$-ideals whose associated forcing is proper may be useful for answering Question \ref{main question analytic} since it indicates how to produce $I^+$ $\borel$ sets. It is should be noted that some restrictions on the type of $\sigma$-ideals considered in Question \ref{main question analytic} is necessary: 

Let $F_{\omega_1}$ denote the countable admissible ordinal equivalence relation defined by $x \ F_{\omega_1} \ y$ if and only if $\omega_1^x = \omega_1^y$. $F_{\omega_1}$ is a thin $\analytic$ equivalence relation with all $\borel$ classes. Thin means that $F_{\omega_1}$ does not have a perfect set of pairwise $F_{\omega_1}$-inequivalent elements. Let $I$ be the $\sigma$-ideal which is $\sigma$-generated by the $F_{\omega_1}$-classes. Suppose there was an $I^+$ $\borel$ set $C$ so that $F_{\omega_1} \upharpoonright C$ is $\borel$. By definition of $I$, each $F_{\omega_1}$-class is in $I$. So since $C$ is $I^+$, $C$ must intersect nontrivially uncountably many classes of $F_{\omega_1}$. So $F_{\omega_1} \upharpoonright C$ has uncountable many classes. Since $F_{\omega_1}$ is thin, there is also no perfect set of $F_{\omega_1} \upharpoonright C$ inequivalent elements. This contradicts Silver's dichotomy (see Fact \ref{silver's dichotomy}).

Of course, $I$ is not proper or even $\omega_1$-preserving: Let $G \subseteq \bbP_I$ be a $\bbP_I$-generic filter over $V$. Fact \ref{name for generic element} implies that $\dot x_\text{gen}[G]$ is not in any ground model coded $\borel$ set in $I$. $\omega_1^{\dot x_\text{gen}[G]}$ can not be a countable admissible ordinal of $V$ since if it was countable then a theorem of Sacks shows that there is a $z \in (\bairespace)^V$ so that $\omega_1^z = \omega_1^{\dot x_\text{gen}[G]}$. Then $x \in [z]_{F_{\omega_1}}$. By definition of $I$, $[z]_{F_{\omega_1}}$ is a $\borel$ set coded in $V$ that belongs to $I$. Hence $\omega_1^{\dot x_\text{gen}[G]}$ must be an uncountable admissible ordinal of $V$, but in $V[G]$, $\omega_1^{\dot x_\text{gen}[G]}$ is a countable admissible ordinal. Hence $\bbP_I$ collapses $\omega_1$.

\Begin{definition}{measure}
A measure $\mu$ on a set $X$ is a nonprincipal ultrafilter on $X$. Nonprincipal means for all $x \in X$, $\{x\} \notin \mu$.

If $\kappa$ is a cardinal, then $\mu$ is $\kappa$-complete if and only if for all $\beta < \kappa$ and sequences $(A_\alpha : \alpha < \beta)$ with each $A_\alpha \in \mu$, $\bigcap_{\alpha < \beta} A_\alpha \in \mu$. $\aleph_1$-completeness is often called countably completeness.

Let $\mathrm{meas}_\kappa(X)$ be the set of all $\kappa$-complete ultrafilter on $X$. 

Suppose $\mu \in \mathrm{meas}_{\aleph_1}({}^{<\omega} X)$. By countably completeness, there is a unique $m$ so that ${}^m X \in \mu$. In this case, $m$ is called the dimension of $\mu$ and this is denoted $\dim(\mu) = m$. 
\end{definition}

\Begin{definition}{tower of measure}
Let $X$ be a set. For $m \leq n < \omega$, let $\pi_{n,m} : {}^n X \rightarrow {}^m X$ be defined by $\pi_{n,m}(f) = f \upharpoonright m$. 

Let $m \leq n < \omega$. Let $\nu$ be a measure of dimension $m$ and $\mu$ be a measure of dimension $n$. $\mu$ is an extension of $\nu$ (or $\nu$ is a projection of $\mu$) if and only if for all $A \in \nu$ with $A \subseteq {}^m X$, $\pi_{n,m}^{-1}[A] \in \mu$. 

A tower of measures over $X$ is a sequence $(\mu_n : n \in \omega)$ so that 

\noindent (i) For all $n$, $\mu_n \in \mathrm{meas}_{\aleph_1}({}^{<\omega} X)$ and $\dim(\mu_n) = n$. 

\noindent (ii) For all $m \leq n < \omega$, $\mu_n$ is an extension of $\mu_m$. 

A tower of measures over $X$, $(\mu_n : n \in \omega)$, is countably complete if and only if for all sequence $(A_n : n \in \omega)$ with the property that for $n \in \omega$, $A_n \in \mu_n$, there exists a $f : \omega \rightarrow X$ so that for all $n \in \omega$, $f \upharpoonright n \in A_n$. 
\end{definition}

\Begin{definition}{trees}
A tree $T$ on $X$ is a subset of ${}^{<\omega} X$ so that if $s \subseteq t$ and $t \in T$, then $s \in T$. 

If $s \in {}^n(X \times Y)$ where $n \in \omega$, then in a natural way, $s$ be may be considered as a pair $(s_0, s_1)$ with $s_0 \in {}^{n}X$ and $s_1 \in {}^nY$. 

Let $T$ be a tree on $X$. The body of $T$, denoted $[T]$, is the set of infinite paths through $T$, that is $[T] = \{f \in {}^\omega X : (\forall n \in \omega)(f\upharpoonright n \in T)\}$. 

Suppose $T$ is a tree on $X \times Y$. For each $s \in {}^{<\omega} X$, define $T^{s} = \{t \in {}^{|s|}Y : (s,t) \in T\}$. If $f \in {}^\omega X$, then define $T^f = \bigcup_{n \in \omega} T^{f\upharpoonright n}$.

Let $T$ be a tree on $X \times Y$, then 
$$p[T] = \{f \in {}^{\omega}X : T^f \text{ is ill-founded} \} = \{f \in {}^{\omega}X : [T^f] \neq \emptyset\}$$
\end{definition}

\Begin{definition}{analytic sets}
For any $k \in \omega$, $A \subseteq {}^k(\bairespace)$ is $\analytic$ if and only if there exists a tree on ${}^k\omega \times \omega$ so that $A = p[T]$. $A \subseteq {}^k(\bairespace)$ is $\coanalytic$ if and only if $A = {}^k(\bairespace) \setminus B$ for some $\analytic$ set $B \subseteq {}^k(\bairespace)$. $A \subseteq {}^k(\bairespace)$ is $\borel$ if and only if $A$ is both $\analytic$ and $\coanalytic$. 
\end{definition}

\Begin{definition}{homogeneous tree}
Let $\gamma$ be an ordinal and $k \in \omega$. A tree $T$ on ${}^k\omega \times \gamma$ is homogeneous if and only if there is a collection $(\mu_s : s \in {}^{<\omega}({}^{k}\omega))$ so that 

\noindent (i) For each $s \in {}^{<\omega}({}^k\omega)$, $\mu_s \in \mathrm{meas}_{\aleph_1}({}^{<\omega}\gamma)$ and concentrates on $T^s$ (that is, $T^s \in \mu_s$).

\noindent (ii) For all $s, t \in {}^{<\omega}({}^k\omega)$, if $s \subseteq t$, then $\mu_t$ is an extension of $\mu_s$.

\noindent (iii) For all $f \in p[T]$, $(u_{f \upharpoonright n} : n \in \omega)$ is a countably complete tower of measures on $\gamma$. 

A collection $(u_s : s \in {}^{<\omega}({}^k\omega))$ which witnesses the homogeneity of $T$ is called a homogeneity system for $T$. 

Let $\kappa$ be a cardinal. The homogeneous tree $T$ is $\kappa$-homogeneous if and only if each $\mu_s$ is $\kappa$-complete.
\end{definition}

\Begin{definition}{homogeneously suslin sets}
For any $k \in \omega$, $A \subseteq {}^k(\bairespace)$ is homogeneously Suslin if and only if there exists an ordinal $\gamma$ and a homogeneous tree on ${}^k\omega \times \gamma$ so that $A = p[T]$.

If the tree $T$ is $\kappa$-homogeneous, then $A$ is said to be $\kappa$-homogeneously Suslin.
\end{definition}

Homogeneously Suslin sets have an important role in the theory of determinacy. In particular, games on $\bairespace$ associated with homogeneously Suslin sets are determined. Later, the homogeneity system of homogeneous trees will be used to show a certain player has a winning strategy in a particular game using techniques that are very similar to the Martin proof of $\analytic$ determinacy from a measurable cardinal.

Below, the basic setting of the relevant games will be described:

\Begin{definition}{games}
Let $X$ be some set. Let $A \subseteq {}^{\omega} X$. The game associated to $A$, denoted $G_A$, is the following: The game has two players, Player 1 and Player 2, who alternatingly take turns playing elements of $X$ with Player 1 playing first. The picture below denotes a partial play where Player 1 plays the sequence $(a_i : i \in \omega)$ and Player 2 plays the sequence $(b_i : i \in \omega)$.
$$\begin{matrix}
a_0 & {} & a_1 & {} & ... & a_{k - 1} & {} \\
\hline 
{} & b_0 & {} & b_1 & ... & {} & b_{k - 1}
\end{matrix}$$
Player 2 is said to win this play of $G_A$ if and only if the infinite sequence $(a_0 b_0 a_1 b_1 ...) \in A$. Otherwise Player 1 wins.

A function $\tau : {}^{<\omega} X \rightarrow X$ is a winning strategy for Player 1 if and only if for all sequence $(b_i : i \in \omega)$ played by Player 2, Player 1 wins by playing $(a_i : i \in \omega)$ where this sequence is defined recursively by $a_{0} = \tau(\emptyset)$ and $a_{k + 1} = \tau(a_0 b_0 ... a_k b_k)$.

A winning stategy $\tau : {}^{<\omega} X \rightarrow X$ for Player 2 is defined similarly.

The game $G_A$ is determined if Player 1 or Player 2 has a winning strategy.
\end{definition}

Let $X$ be a set. ${}^{\omega} X$ is given the topology with basis $\{U_s : s \in {}^{<\omega}X\}$, where $U_s = \{f \in {}^{\omega}X : s \subseteq f\}$. 

\Begin{fact}{open determinacy}
(\cite{Infinite-Games-With-Perfect-Information}, Gale-Stewart) If $A \subseteq {}^\omega X$ is open, then $G_A$ is determined. Hence if $A$ is closed, then $G_A$ is also determined.
\end{fact}

\section{The Game}\label{The Game}
\Begin{definition}{section of relation}
Let $R$ be a relation on $(\bairespace)^2$. Let $R_x = \{y : (x,y) \in R\}$ and $R^x = \{y : (y,x) \in R\}$
\end{definition}

The following results will be stated using the verticle sections $R_x$; however, the results hold using horizontal sections with the appropriate changes.

\Begin{definition}{homogeneously suslin relation}
Let $S$ be a homogeneous tree on $\omega \times \omega \times \gamma$, where $\gamma$ is some ordinal.

Denote $p[S]$ by $R_S$. Denote $(\bairespace \times \bairespace) \setminus p[S] = R^S$. 
\end{definition}

\Begin{definition}{Assumption A}
Let $S$ be a homogeneously tree on $\omega \times \omega \times \gamma$ for some ordinal $\gamma$. Let $I$ be a $\sigma$-ideal on $\bairespace$ so that $\bbP_I$ is proper.

Assumption $\mathsf{A_\Sigma}$ asserts that $1_{\bbP_I} \forces_{\bbP_I} \check S$ is a homogeneous tree.

Assumption $\mathsf{A_\Pi}$ asserts $1_{\bbP_I} \forces_{\bbP_I} \check S$ is a homogeneous tree.
\end{definition}

Assumption $\mathsf{A_\Sigma}$ and $\mathsf{A_\Pi}$ just asserts that the tree $S$ remains homogeneous in $\bbP_I$-generic extensions. (Since the completeness of countably complete measures is a measurable cardinal and $|\bbP_I|$ is always less than a measurable cardinal under $\mathsf{AC}$, this is always true under $\mathsf{AC}$.) 

\Begin{definition}{D formula}
Let $S$ be a homogeneously Suslin tree on $\omega \times \omega \times \gamma$ for some ordinal $\gamma$. Let $I$ be a $\sigma$-ideal on $\bairespace$ such that $\bbP_I$ is a proper forcing.

Let $D_\Sigma$ be the formula on $\bairespace \times \bairespace$ asserting: 
$$D_\Sigma(x,T) \Leftrightarrow (T \text{ is tree on } \omega \times \omega) \wedge (\forall y)(R_S(x,y) \Leftrightarrow T^y \text{ is ill-founded})$$

Let $D_\Pi$ be the formula on $\bairespace \times \bairespace$ asserting:
$$D_\Pi(x,T) \Leftrightarrow (T \text{ is a tree on } \omega \times \omega) \wedge (\forall y)(\neg(R^S(x,y)) \Leftrightarrow T^y \text{ is ill-founded})$$
\end{definition}

If $D_\Sigma(x,T)$ holds, then $T$ is a tree which witnesses $(R_S)_x$ is $\analytic$. Similarly, if $D_\Pi(x,T)$ holds, then $T$ is a tree which witnesses $\bairespace \setminus (R^S)_x$ is $\analytic$, i.e. $(R^S)_x$ is $\coanalytic$.

\Begin{definition}{Assumption B}
Let $S$ be a homogeneously Suslin tree on $\omega \times \omega \times \gamma$ for some ordinal $\gamma$. Let $I$ be a $\sigma$-ideal on $\bairespace$ such that $\bbP_I$ is a proper forcing.

Let assumption $\mathsf{B_\Sigma}$ say: $(\forall x)(\exists T)D_\Sigma(x,T)$ and $1_{\bbP_I} \forces_{\bbP_I} (\forall x)(\exists T)D_\Sigma(x,T)$.

Let assumption $\mathsf{B_\Pi}$ say: $(\forall x)(\exists T)D_\Pi(x,T)$ and $1_{\bbP_I} \forces_{\bbP_I} (\forall x)(\exists T)D_\Pi(x,T)$.
\end{definition}

Assumption $\mathsf{B_\Sigma}$ states that all $R_S$ sections are $\analytic$ and all $R_S$ sections remain $\analytic$ in $\bbP_I$-generic extensions. Similarly, assumption $\mathsf{B_\Pi}$ states that all $R^S$ sections are $\coanalytic$ and all $R^S$ sections remain $\coanalytic$ in $\bbP_I$-generic extensions.

\Begin{definition}{Assumption C}
Let $S$ be a homogeneously Suslin tree on $\omega \times \omega \times \gamma$ for some ordinal $\gamma$. Let $I$ be a $\sigma$-ideal on $\bairespace$ such that $\bbP_I$ is a proper forcing. 

Let assumption $\mathsf{C_\Sigma}$ state: There is an ordinal $\epsilon$ and a tree $U$ on $\omega \times \omega \times \epsilon$ so that $p[U] = \{(x,T) : D_\Sigma(x,T)\}$ and $1_{\bbP_I} \forces_{\bbP_I} p[\check U] = \{(x,T) : D_\Sigma(x,T)\}$. 

Let assumption $\mathsf{C_\Pi}$ state: There is an ordinal $\epsilon$ and a tree $U$ on $\omega \times \omega \times \epsilon$ so that $p[U] = \{(x,T) : D_\Pi(x,T)\}$ and $1_{\bbP_I} \forces_{\bbP_I} p[\check U] = \{(x,T) : D_\Pi(x,T)\}$. 
\end{definition}

Assumption $\mathsf{C_\Sigma}$ states that the set defined by $D_\Sigma$ has a tree representation that continues to represent the formula $D_\Sigma$ in $\bbP_I$-generic extensions. $\mathsf{C_\Pi}$ is similar.

The following shows under certain assumptions a more general canonicalization property holds for relations. \cite{Borel-Canonization-of-Analytic-Sets} defines this phenomenon as the rectangular canonization property.

\Begin{theorem}{homogeneous suslin equiv with all analytic class analytic somewhere}
Let $\gamma$ be an ordinal. Let $S$ be a homogeneous tree on $\omega \times \omega \times \gamma$. Let $I$ be a $\sigma$-ideal on $\bairespace$ so that $\bbP_I$ is proper. Assume $\mathsf{A_\Sigma}$, $\mathsf{B_\Sigma}$, and $\mathsf{C_\Sigma}$ hold for $S$ and $I$. 

Then for any $I^+$ $\borel$ set $B \subseteq \bairespace$, there exists an $I^+$ $\borel$ set $C \subseteq B$ so that $R_S \cap (C \times \bairespace)$ is an $\analytic$ relation.
\end{theorem}

\begin{proof}
Let $U$ be the tree on $\omega \times \omega \times \epsilon$ witnessing $\mathsf{C_\Sigma}$ for $S$ and $I$. 

Let $M \prec H_\Xi$ be a countable elementary substructure with $\Xi$ sufficiently large and $B,I,\bbP_I, S, U \in M$. 
\\*
\\*\underline{Claim 1} : Let $g$ be $\bbP_I$-generic over $M$. If $x,T \in M[g]$ and $M[g] \models D_\Sigma(x,T)$, then $V \models D_\Sigma(x,T)$. 

Proof of Claim 1: By assumption $C_\Sigma$ for $S$ and $I$ and the fact that $M \prec H_\Xi$, $M[g] \models D_\Sigma(x,T)$ implies $M[g] \models (x,T) \in p[U]$. There exists some $f \in M[g]$ with $f : \omega \rightarrow \epsilon$ so that $M[g] \models (x,T,f) \in [U]$. Hence for each $n \in \omega$, $M[g] \models (x \upharpoonright n, T \upharpoonright n, f \upharpoonright n) \in U$. For each $n \in \omega$, $(x \upharpoonright n, T \upharpoonright n, f \upharpoonright n) \in M$. So by absoluteness, $M \models (x \upharpoonright n, T \upharpoonright n, f \upharpoonright n) \in U$. For all $n \in \omega$, $V \models (x \upharpoonright n, T \upharpoonright n, f \upharpoonright n) \in U$. $V \models (x,T) \in p[U]$. $V \models D_\Sigma(x,T)$. 
\\*
\\*\indent Now fix a $g \in \bairespace$ so that $g$ is $\bbP_I$-generic over $M$.

As $M \prec H_\Xi$, $M \models (\forall x)(\exists T)D_\Sigma(x,T)$. $M[g] \models (\forall x)(\exists T)D_\Sigma(x,T)$ by assumption $B_\Sigma$ and the fact that $M \prec H_\Xi$. So fix a tree $T$ on $\omega \times \omega$ so that $M[g] \models D_\Sigma(g,T)$. 

Consider the following game $G^{g,T}$: 
$$\begin{matrix}
m_0,n_0 & {} & m_1,n_1 & {} & ... & m_{k - 1}, n_{k - 1} & {} \\
\hline 
{} & \alpha_0 & {} & \alpha_1 & ... & {} & \alpha_{k - 1}
\end{matrix}$$
The rules are:

\noindent (1) Player 1 plays $m_i, n_i \in \omega$. Player 2 plays $\alpha_i < \gamma$.

\noindent (2) $(m_0 ... m_{k - 1}, n_0 ... n_{k - 1}) \in T$

\noindent (3) $(g \upharpoonright k, m_0 ... m_{k - 1}, \alpha_0 ... \alpha_{k - 1}) \in S$. 

The first player to violate these rules loses. If the game continues forever, then Player 2 wins. 
\\*
\\* \underline{Claim 2} : In $M[g]$, Player 2 has a winning strategy in the game $G^{g,T}$. 

Proof of Claim 2: By an appropriate coding, $G^{g,T}$ is equivalent to a game $G_A$, where $A \subseteq {}^{\omega} \gamma$ is a closed subset. 

Suppose Player 2 does not have a winning strategy. By closed determinacy (Fact \ref{open determinacy}), Player 1 must have a winning strategy $\tau^*$. 

By assumption $\mathsf{A_\Sigma}$, $S$ is a homogeneous tree in $M[g]$. Let $(\mu_t : t \in {}^{<\omega}(\omega \times \omega))$ be a homogeneity system witnessing the homogeneity of $S$. 

Now two sequences of natural numbers, $(a_i : i \in \omega)$ and $(b_i : i \in \omega)$, and a sequence $(A_n : n \in \omega)$ so that $A_n \subseteq {}^{n}\gamma$ will be constructed by recursion:

Let $a_0,b_0 \in \omega$ so that $(a_0,b_0) = \tau^*(\emptyset)$. Let $A_0 = \{\emptyset\}$.

Suppose $a_0, ..., a_{k-1}$, $b_0, ..., b_{k - 1}$, and $A_0, ..., A_{k -1}$ has been constructed. Define the function 
$$h_{k} : S^{(g \upharpoonright k, a_0...a_{k - 1})} \rightarrow \omega \times \omega$$
defined by
$$h_{k}(\beta_0...\beta_{k - 1}) = \tau^*(a_0, b_0, \beta_0,...,a_{k - 1}, b_{k - 1}, \beta_{k - 1})$$
$\mu_{(g\upharpoonright k,a_0...a_{k - 1})}$ concentrates on $S^{(g\upharpoonright k, a_0...a_{k - 1})}$ and is countably complete; therefore, there is a unique $(a_k, b_k)$ so that 
$$h^{-1}_{k}[\{(a_k,b_k)\}] \in \mu_{(g\upharpoonright k, a_0...a_{k - 1})}$$
Let $A_{k} = h^{-1}_{k}[\{(a_k,b_k)\}]$. 

This completes the construction of $(a_i : i \in \omega)$, $(b_i : i \in \omega)$, and $(A_i : i \in \omega)$.

Let $L \in {}^\omega (\omega \times \omega)$ be such that for all $i \in \omega$, $L(i) = (a_i,b_i)$. Note that $L \in [T]$. To see this, suppose not. Then there is some least $k \in \omega$ so that $L \upharpoonright (k + 1) = (a_0...a_k, b_0...b_k) \notin T$. For $i \leq k$, define $\mu_i = \mu_{g \upharpoonright i, a_0...a_{i - 1}}$. For $0 \leq i \leq j \leq k$, let $\pi_{j,i} : {}^{j}\gamma \rightarrow {}^i\gamma$ be defined by $\pi_{j,i}(s) = s \upharpoonright i$. By definition of the homogeneity system for $S$, for $0 \leq i \leq j \leq k$, $\mu_j$ is an extension of $\mu_i$. Hence for all $0 \leq i \leq k$, $\pi_{k,i}^{-1}[A_i] \in \mu_k$. By countable completeness of $\mu_k$, $\bigcap_{0 \leq i \leq k} \pi_{k,i}^{-1}[A_i] \in \mu_k$. Let $(\beta_0...\beta_{k-1}) \in \bigcap_{0 \leq i \leq k} \pi_{k,i}^{-1}[A_i]$. Consider the following play of $G^{g,T}$ where player 1 uses the strategy $\tau^*$ and Player 2 plays $(\beta_0...\beta_{k - 1})$: 
$$\begin{matrix}
a_0,b_0 & {} & a_1,b_1 & {} & ... & a_{k - 1}, b_{k - 1} & {} & a_k,b_k \\
\hline 
{} & \beta_0 & {} & \beta_1 & ... & {} & \beta_{k - 1} & {}
\end{matrix}$$
Note that for all $0 \leq i \leq k$,  $(\beta_0...\beta_{i - 1}) \in A_i = h_i^{-1}[\{(a_i,b_i)\}] \subseteq S^{(g\upharpoonright i, a_0...a_{i - 1})}$. So rule (3) of the game $G^{g,T}$ is not violated by Player 2. However, $(a_0...a_k, b_0...b_k) = L \upharpoonright (k + 1) \notin T$. Player 1 violates rule (2) and is the first player to violate any rules. Player 1 loses this game. This contradicts the assumption that $\tau^*$ is a winning strategy for Player 1. So this completes the proof that $L \in [T]$. 

Let $\mathfrak{a} = (a_i : i \in \omega)$. Since $L \in [T]$ and $D_\Sigma(g,T)$, this implies that $R_S(g,\mathfrak{a})$. 

Now let $J \in {}^{\omega} (\omega \times \omega)$ be such that for all $k \in \omega$, $J \upharpoonright k = (g \upharpoonright k, a_0 ... a_{k - 1})$. Then by definition of $S$, $J \in p[S]$. Since $S$ is a homogeneous tree via $(u_t : t \in {}^{<\omega}(\omega \times \omega))$, $(\mu_{J\upharpoonright k} : k \in \omega)$ is a countably complete tower of measures. 

Each $A_k \in \mu_{g\upharpoonright k, a_0...a_{k -1}} = \mu_{J \upharpoonright k}$. So by the countable completeness of the tower, there exists some $\Phi : \omega \rightarrow \gamma$ so that for all $k \in \omega$, $\Phi \upharpoonright k \in A_k$. Now consider the play of $G^{g,T}$ where Player 1 uses its winning strategy $\tau^*$ and Player 2 plays $\Phi$. By construction of the sequences $(a_i : i \in \omega)$, $(b_i, i \in \omega)$, and $(A_i : i \in \omega)$, the game looks as follows:
$$\begin{matrix}
a_0,b_0 & {} & a_0,b_1 & {} & ... & a_{k - 1}, b_{k - 1} & {} \\
\hline 
{} & \Phi(0) & {} & \Phi(1) & ... & {} & \Phi(k - 1)
\end{matrix}$$
Neither players violate any rules in this play. Hence the game continues forever, and so Player 2 wins this play of $G^{g,T}$. This contradicts the fact that $\tau^*$ was a winning strategy for Player 1.

So Player 1 could not have had a winning strategy. Player 2 must have a winning strategy in $G^{g,T}$. This completes the proof of Claim 2.
\\*
\\*\indent By Claim 2, Player 2 has a winning strategy $\tau \in M[g]$.
\\*
\\* \underline{Claim 3} : $\tau$ is a winning strategy for $G^{g,T}$ in $V$. 

Proof of Claim 3: Suppose the following is a play of $G^{g,T}$ in which Player 2 uses $\tau$ and loses
$$\begin{matrix}
m_0,n_0 & {} & m_1,n_1 & {} & ... & m_{k - 1}, n_{k - 1} & {} \\
\hline 
{} & \alpha_0 & {} & \alpha_1 & ... & {} & \alpha_{k - 1}
\end{matrix}$$
Since $\tau \in M[g]$ and ${}^{<\omega}\omega \subseteq M[g]$, this entire finite play belongs to $M[g]$. So, Player 2 loses this game in $M[g]$, as well. This contradicts $\tau$ being a winning strategy in $M[g]$. This completes the proof of Claim 3.
\\*
\\*\indent \underline{Claim 4} : For all $y \in \bairespace$, $R_S(g,y)$ if and only if $(S \cap M)^{(g,y)}$ is ill-founded.

Proof of Claim 4: By Claim 1, $M[g] \models D_\Sigma(g,T)$ implies $V \models D_\Sigma(g,T)$. Hence in $V$, $T$ gives the $\analytic$ definition of $(R_S)_g$. 

Suppose $R_S(g,y)$. Then $T^y$ is ill-founded. Let $f \in [T^y]$. Consider the following play of the game $G^{g,T}$ where Player 1 plays $y$ and $f$, and Player 2 responds using its winning strategy $\tau$. 
$$\begin{matrix}
y(0),f(0) & {} & y(1),f(1) & {} & ... & y(k - 1), f(k - 1) & {} \\
\hline 
{} & \alpha_0 & {} & \alpha_1 & ... & {} & \alpha_{k - 1}
\end{matrix}$$
Since $f \in [T^y]$, Player 1 can not lose. Since $\tau$ is a winning strategy for Player 2, Player 2 also does not lose at a finite stage. Hence Player 2 wins by having the game continue forever. Let $\Phi : \omega \rightarrow \gamma$ be the sequence coming from Player 2's response, i.e. for all $k$, $\Phi(k) = \alpha_k$. 

Since $\tau \in M[g]$ and $\finNaturalSequence \subseteq M[g]$, each finite partial play of $G^{g,T}$ above belongs to $M[g]$. Hence $\Phi \upharpoonright k \in M[g]$ for all $k \in \omega$. As $\text{On}^M = \text{On}^{M[g]}$ because $\bbP_I$ is proper, $(g \upharpoonright k, y \upharpoonright k, \Phi \upharpoonright k) \in (S \cap M)$ for all $k \in \omega$.

It has been shown that $R_S(g,y)$ implies $(S \cap M)^{(g,y)}$ is ill-founded. 

Of course, if $(S \cap M)^{(g,y)}$ is ill-founded, then $S^{(g,y)}$ is ill-founded. By definition, $R_S(g,y)$. 

This completes the proof of Claim 4.
\\*
\\*\indent Let $b : \omega \rightarrow \text{On}^M$ be a bijection. Define a new tree $S'$ on $\omega \times \omega \times \omega$ by $(s_1,s_2,s_3) \in S' \Leftrightarrow (s_1,s_2, b \circ s_3) \in S$. 

By Fact \ref{properness equivalence}, let $C \subseteq B$ be the $I^+$ $\borel$ set of $\bbP_I$-generic reals over $M$ inside $B$. By Claim 4, for all $y \in \bairespace$, $R_S(g,y) \Leftrightarrow (S')^{(g,y)}$ is ill-founded. $R_S \cap (C \times \bairespace)$ is $\analytic$. The proof of the theorem is complete.
\end{proof}

\Begin{theorem}{homogeneous suslin equiv with all coanalytic class coanalytic somewhere}
Let $\gamma$ be an ordinal. Let $S$ be a homogeneous tree on $\omega \times \omega \times \gamma$. Let $I$ be a $\sigma$-ideal on $\bairespace$ so that $\bbP_I$ is proper. Assume $\mathsf{A_\Pi}$, $\mathsf{B_\Pi}$, and $\mathsf{C_\Pi}$ holds for $S$ and $I$. 

Then for any $I^+$ $\borel$ set $B \subseteq \bairespace$, there exists an $I^+$ $\borel$ set $C \subseteq B$ so that $R^S \cap (C \times \bairespace)$ is a $\coanalytic$ relation.
\end{theorem}

\begin{proof}
The proof of this is very similar to the proof of Theorem \ref{homogeneous suslin equiv with all analytic class analytic somewhere}.
\end{proof}

\Begin{theorem}{homogeneous suslin equiv with all borel classes borel somewhere}
Let $\gamma$ and $\nu$ be ordinals. Let $S$ be a homogeneous tree on $\omega \times \omega \times \gamma$. Let $U$ be a homogeneous tree on $\omega \times \omega \times \nu$. Suppose $p[S] = (\bairespace \times \bairespace) \setminus p[U]$. Let $R = R_S = R^U$. Let $I$ be a $\sigma$-ideal on $\bairespace$ such that $\bbP_I$ is a proper forcing. Suppose $\mathsf{A_\Sigma}$, $\mathsf{B_\Sigma}$, and $\mathsf{C_\Sigma}$ holds for $S$ and $I$. Suppose $\mathsf{A_\Pi}$, $\mathsf{B_\Pi}$, and $\mathsf{C_\Pi}$ holds for $U$ and $I$. 

Then for any $I^+$ $\borel$ set $B \subseteq \bairespace$, there exists an $I^+$ $\borel$ set $C \subseteq B$ so that $R \cap (C \times \bairespace)$ is a $\borel$ relation.
\end{theorem}

\begin{proof}
By Theorem \ref{homogeneous suslin equiv with all analytic class analytic somewhere}, there is some $I^+$ $\borel$ set $C' \subseteq B$ so that $R \cap (C' \times \bairespace)$ is $\analytic$. By Theorem \ref{homogeneous suslin equiv with all coanalytic class coanalytic somewhere}, there is some $I^+$ $\borel$ set $C \subseteq C'$ so that $R \cap (C \times \bairespace)$ is $\coanalytic$. Therefore, $R \cap (C \times \bairespace)$ is $\borel$. 
\end{proof}

If the above assumptions holds and $R_S = E$ defines an equivalence relation with all $\analytic$ classes, then there is some $I^+$ $\borel$ set so that $E \upharpoonright C$ is an $\analytic$ equivalence relation.

Simialarly, suppose $R_S = G$ is a graph on $\bairespace$. Then $G_x = \{y : x \ G \ y\}$ is the set of neighbors of $x$. Suppose $G_x$ is $\analytic$ for all $x$. Then there is an $I^+$ $\borel$ set $C$ so that the induced subgraph $G \upharpoonright C$ is an $\analytic$ graph.

Since equivalence relations were the original motivation, the rest of the paper will focus on equivalence relations; however, all the results holds for graphs and relations with the appropriate sections. 

\section{Canonicalization for Equivalence Relations in $L(\bbR)$}\label{Canonicalization}
This section will provide a brief description of the theory of tree representations of subsets of $\bairespace$ and absoluteness. This will be used to indicate some circumstances in which the assumptions $\mathsf{A_\Sigma}$, $\mathsf{B_\Sigma}$, $\mathsf{C_\Sigma}$, $\mathsf{A_\Pi}$, $\mathsf{B_\Pi}$, and $\mathsf{C_\Pi}$ hold. The results of the previous section will be applied to some familiar classes of equivalence relations. The following discussion is in $\mathsf{ZF} + \mathsf{DC}$ until it is explicitly mentioned that $\mathsf{AC}$ will be assumed.

\Begin{definition}{weakly homogeneous systems}
Let $\kappa$ be a cardinal. A $\kappa$-weak homogeneity system with support some ordinal $\gamma$ is a sequence of $\kappa$-complete measures on ${}^{<\omega} \gamma$, $\bar{\mu} = (\mu_s : s \in \finNaturalSequence)$, so that 

(i) If $s \neq t$, then $\mu_s \neq \mu_t$. 

(ii) $\dim(\mu_s) \leq |s|$.

(iii) If $\mu_s$ is an extension of some measure $\nu$, then there exists some $k < |s|$ so that $\mu_{s \upharpoonright k} = \nu$. 

Define $W_{\bar\mu}$ by
$$W_{\bar\mu} = \{x \in \bairespace : (\exists f \in \bairespace)(f \text{ is an increasing  sequence} \wedge (\mu_{x \upharpoonright f(k)} : k \in \omega) \text{ is a countably complete tower})\}$$

A set $A \subseteq \bairespace$ is $\kappa$-weakly homogeneous if and only there is a $\kappa$-weak homogenity system $\bar\mu$ so that $A = W_{\bar\mu}$. 
\end{definition}

\Begin{definition}{weakly homogeneous trees}
Let $\gamma$ be an ordinal. A tree on $\omega \times \gamma$ is $\kappa$-weakly homogeneous if and only there is some $\kappa$-weak homogeneity system $\bar\mu = (\mu_s : s \in \finNaturalSequence)$ so that $p[T] = W_{\bar\mu}$ and for all $s \in \finNaturalSequence$, there is some $k \leq |s|$ so that $\mu_s$ concentrates on $T^{s \upharpoonright k}$. 

$A \subseteq \bairespace$ is $\kappa$-weakly homogeneously Suslin if and only if $A = p[T]$ for some tree $T$ which is $\kappa$-weakly homogeneous.
\end{definition}

\Begin{fact}{weakly homogeneous and weakly homogeneously suslin}
If $\bar\mu = (\mu_s : s \in \finNaturalSequence)$ is a $\kappa$-weak homogeneity system with support $\gamma$, then there is a tree $T$ on $\omega \times \gamma$ so that $\bar\mu$ witnesses $T$ is $\kappa$-weakly homogeneously Suslin.

Hence a set is $\kappa$-weakly homogeneous if and only if it is $\kappa$-weakly homogeneously Suslin.
\end{fact}

\begin{proof}
See \cite{The-Derived-Model-Theorem}, Proposition 1.12. 
\end{proof}

\Begin{definition}{homogeneous system of elem embedding}
Let $\mu$ be a countably complete measure on ${}^{<\omega} X$. Let $M_\mu$ be the Mostowski collapse of the the ultrapower $\mathrm{Ult}(V,\mu)$. Let $j_\mu : V \rightarrow M_\mu$ be the composition of the ultrapower map and the Mostowski collapse map.

Suppose $\nu$ and $\mu$ are countably complete measures on ${}^{<\omega} X$. Suppose for some $m \leq n$, $\dim(\mu) = m$ and $\dim(\nu) = n$, and $\nu$ is an extension of $\mu$. Define $\Lambda_{m,n} : {}^{{}^m X} V \rightarrow {}^{{}^nX} V$ by $\Lambda_{m,n}(f)(s) = f(s \upharpoonright m)$ for each $s \in {}^{n}X$. Define an elementary embedding $\mathrm{Ult}(V,\nu) \rightarrow \mathrm{Ult}(V,\mu)$ by $[f]_\nu \mapsto [\Lambda_{m,n}(f)]_\mu$. This induces an elementary embedding $j_{\nu,\mu} : M_\nu \rightarrow M_\mu$. 
\end{definition}

\Begin{definition}{martin-solovay tree}
Let $\gamma$ and $\theta$ be ordinals. Let $\bar\mu = (\mu_s : s \in \finNaturalSequence)$ be a weak homogeneity system with support $\gamma$. The Martin-Solovay tree with respect to $\bar\mu$ below $\theta$, denoted $\mathrm{MS}_\theta(\bar\mu)$, is a tree on $\omega \times \theta$ defined by: for all $s \in \finNaturalSequence$ and $h \in {}^{|s|}\theta$
$$(s,h) \in \mathrm{MS}_\theta(\bar\mu) \Leftrightarrow (\forall i < j < |s|)(\mu_{s \upharpoonright j} \text{ is an extension of } \mu_{s \upharpoonright i} \Rightarrow j_{\mu_{s \upharpoonright i}, \mu_{s \upharpoonright j}}(h(i)) > h(j))$$
\end{definition}

If $(u_n : n \in \omega)$ is a tower of measure, then the tower is countably complete if and only if the directed limit of the directed system $(M_{\mu_i} : j_{\mu_i,\mu_j} : i < j < \omega)$ is well-founded. Suppose $x \in p[\mathrm{MS}_\theta(\bar\mu)]$. If $(x,\Phi) \in [\mathrm{MS}_\theta(\bar\mu)]$, then $\Phi$ witnesses in a continuous way that the directed limit model is ill-founded. This shows that $x \in p[\mathrm{MS}_\theta(\bar\mu)]$ implies that $x \notin W_{\bar\mu}$. In fact, the converse is also true giving the following result:

\Begin{fact}{martin-solovay projection}
($\mathsf{ZF + DC}$) Let $\kappa$ be a cardinal. Suppose $\bar\mu$ is a $\kappa$-weak homogeneity system with support $\gamma$. Then if $\theta > |\gamma|^+$, then $p[\mathrm{MS}_\theta(\bar\mu)] = \bairespace \setminus W_{\bar\mu}$. 
\end{fact}

\begin{proof}
See \cite{The-Derived-Model-Theorem} Lemma 1.19, \cite{Stationary-Tower} Fact 1.3.12, or \cite{Suslin-Cardinals-Partition-Properties} Theorem 4.10.
\end{proof}

Let $\mu$ be a $\kappa$-complete ultrafilter on some set $X$. Let $\bbP$ be a forcing with $|\bbP| < \kappa$. Let $G \subseteq \bbP$ be $\bbP$-generic over $V$. It can be shown that if $f^* : X \rightarrow V$ is a function in $V[G]$, then there is a function $f \in V$ and $A \in \mu$ so that $V[G] \models (\forall x \in A)(f(x) = f^*(x))$. 

In $V[G]$, define $\mu^* \subseteq \mathcal{P}(X)$ by $A \in \mu^*$ if and only there exists a $B \in \mu$ so that $B \subseteq A$. In $V[G]$, $\mu^*$ is a $\kappa$-complete ultrafilter on $X$. Let $M_{\mu^*}$ denote the Mostowski collapse of $\mathrm{Ult}(V[G], \mu^*)$. Let $j^*_{\mu^*} : V[G] \rightarrow M_{\mu^*}$ be the induced elementary embedding.

In $V[G]$, $\mathrm{Ult}(V,\mu)$ can be embedded into $\mathrm{Ult}(V[G], \mu^*)$ as follows: for all $f \in ({}^XV) \cap V$, $[f]_\mu \mapsto [f]_{\mu^*}$. If $f \in ({}^X V) \cap V$ and $g' \in {}^X V[G]$ are such that $\mathrm{Ult}(V[G], \mu^*) \models [g']_{\mu^*} \in [f]_{\mu^*}$, then $\{x \in X : g'(x) \in f(x)\} \in \mu^*$. Therefore, one can find a $g^* \in V[G]$ so that $g^* : X \rightarrow V$ and $[g']_{\mu^*} = [g^*]_{\mu^*}$. By the above observation, one can find a $g \in V$ so that $[g]_{\mu^*} = [g^*]_{\mu^*} = [g']_{\mu^*}$. This shows that $\mathrm{Ult}(V,\mu)$ is identified (via the embedding above) as an $\in^{\mathrm{Ult}(V[G], \mu^*)}$-initial segment of $\mathrm{Ult}(V[G], \mu^*)$. After Mostowski collapsing the ultrapowers, it can be seen that $j_{\mu^*} \upharpoonright M_\mu = j_{\mu}$. 

Suppose $\bar\mu = (\mu_s : s \in \finNaturalSequence)$ is a $\kappa$-weak homogeneity system. Denote $\bar\mu^* = (\mu_s^* : s \in \finNaturalSequence)$. $\bar\mu^*$ is a $\kappa$-weak homogeneity system. From the construction, the Martin-Solovay trees depends only on $j_{\mu_s^*} \upharpoonright \text{ON}$. So by the above discussion, $\mathrm{MS}_\theta(\bar\mu)^V = \mathrm{MS}_\theta(\bar\mu^*)^{V[G]}$. Hence Fact \ref{martin-solovay projection} implies that $V[G] \models p[\mathrm{MS}_\theta(\bar{u})] = \bairespace \setminus W_{\bar\mu^*}$. 

(The above argument can be applied to a $\kappa$-homogeneous tree $S$ and its witnessing $\kappa$-homogeneity system $\bar{\mu}$ to show that if $|\bbP| < \kappa$, then $\bar{\mu}^*$ is a $\kappa$-weak homogeneity system for $S$ in $V[G]$. Assuming the axiom of choice, this shows assumption $\mathsf{A_\Sigma}$ and $\mathsf{A_\Pi}$.)

Now suppose that $T$ is a $\kappa$-weakly homogeneous tree on $\omega \times \alpha$ witnessed by the $\kappa$-weak homogeneity system $\bar\mu$. This gives that $p[T] = W_{\bar\mu}$. One seeks to show that $p[\mathrm{MS}_\theta(\bar\mu)^V]$ continues to represent $\bairespace \setminus p[T]$ in $V[G]$. By the previous paragraph, it suffices to show that $V[G] \models p[T] = W_{\bar\mu^*}$: If $x \in W_{\bar\mu^*}$, then there is an increasing function $f : \omega \rightarrow \omega$ so that $(\mu_{x \upharpoonright f(n)}^* : n \in \omega)$ is a countably complete tower. For all $n$, $\mu_{x \upharpoonright f(n)}^*$ concentrates on $T^{x \upharpoonright n}$. So by countably completeness, there is a path $\Phi \in [T^x]$. So $x \in p[T]$. Conversely, suppose $x \in p[T]$. Fact \ref{martin-solovay projection} implies that in $V$, $T$ and $\mathrm{MS}_\theta(\bar\mu)$ are complementing trees. By the absoluteness of well-foundedness, $V[G] \models \emptyset= p[T] \cap p[\mathrm{MS}_\theta(\bar\mu)] = p[T] \cap p[\mathrm{MS}_\theta(\bar\mu^*)]$. So $x \notin p[\mathrm{MS}_\theta(\bar{\mu}^*)]$. Then applying Fact \ref{martin-solovay projection} in $V[G]$ to the weak homogeneity system $\bar\mu^*$, one obtains that $x \in W_{\bar\mu^*}$.

So in summary:

\Begin{fact}{martin-solovay generically correct represention of complement}
($\mathsf{ZF + DC}$) Let $\kappa$ be a cardinal. Let $T$ be a $\kappa$-weakly homogeneous tree on $\omega \times \gamma$, for some ordinal $\gamma$, with $\kappa$-weak homogeneity system $\bar{\mu}$. Let $\theta > |\gamma|^+$. Let $\bbP$ be a forcing with $|\bbP| < \kappa$ and $G \subseteq \bbP$ be $\bbP$-generic over $V$.

$V[G] \models \mathrm{MS}_\theta(\bar\mu^*) = \mathrm{MS}_\theta(\bar\mu)^V$. $V[G] \models p[\mathrm{MS}_\theta(\bar\mu)^V] = \bairespace \setminus p[T]$. 
\end{fact}

\begin{proof}
See \cite{The-Derived-Model-Theorem}, Section 1 and especially Lemma 1.19. Also see \cite{Stationary-Tower}, Section 1.3.
\end{proof}

So if $T$ is $\kappa$-weakly homogeneous, an appropriate Martin-Solovay tree will continue to represent the complement of $p[T]$ in generic extensions by forcings of cardinal less than $\kappa$. The Martin-Solovay trees give the generically-correct tree representations for complements of $\kappa$-weakly homogeneously sets. However, the formulas $D_\Sigma$ and $D_\Pi$ involve more negations and quantifications over $\bairespace$. Multiple iterations of the Martin-Solovay construction will be needed. The following results are useful for continuing the Martin-Solovay construction of generically-correct tree representation for more complex sets. In addition, these results will also imply that these representations are also homogeneously Suslin. Until the end of this section, the axiom of choice will be assumed.

\Begin{definition}{quantifier notation}
If $B \subseteq {}^k(\bairespace) \times \bairespace$. Denote 
$$\exists^\R B = \{x : (\exists y)((x,y) \in B)\}$$
$$\forall^\R B = \{x : (\forall y)((x,y) \in B)\}$$
If $A \subseteq {}^k(\bairespace)$, then denote
$$\neg A = {}^k(\bairespace) \setminus A$$
\end{definition}

\Begin{fact}{weakly homogeneous and homogeneous sets}
Let $A \subseteq \bairespace$. $A$ is $\kappa$-weakly homogeneously Suslin if and only if there is a $\kappa$-homogeneously Suslin set $B \subseteq \bairespace \times \bairespace$ so that $A = \exists^\R B$. 
\end{fact}

\begin{proof}
See \cite{The-Derived-Model-Theorem}, Proposition 1.10.
\end{proof}

A Woodin cardinal is a technical large cardinal which has been very useful in descriptive set theory. (See \cite{Stationary-Tower}, Section 1.5 for more information about Woodin cardinals.)

\Begin{fact}{martin-steel propogation}
Let $\delta$ be a Woodin cardinal. Let $\bar\mu = (\mu_s : s \in \finNaturalSequence)$ be a $\delta^+$-weak homogeneity system with support $\gamma \in \mathrm{ON}$. Then for sufficiently large $\theta$, $\mathrm{MS}_\theta(\bar\mu)$ is $\kappa$-homogeneous for all $\kappa < \delta$.
\end{fact}

\begin{proof}
See \cite{Proof-of-Projective-Determinacy}.
\end{proof}

\Begin{definition}{homogeneous class notation}
If $\kappa$ is a cardinal, then let $\mathrm{Hom}_\kappa$ be the collection of $\kappa$-homogenously Suslin subsets of $\bairespace$. Let $\mathrm{Hom}_{<\kappa} = \bigcap_{\gamma < \kappa} \mathrm{Hom}_\gamma$. 
\end{definition}

The following are some well-known results on what sets can be in $\mathrm{Hom}_{<\lambda}$ when $\lambda$ is limit of Woodin cardinals. 

\Begin{fact}{closure of Hom}
(Martin-Steel) Let $\lambda$ be a limit of Woodin cardinals. Then $\mathrm{Hom}_{<\lambda}$ is closed under complements and $\forall^\bbR$. 
\end{fact}

\begin{proof}
Let $A \in \mathrm{Hom}_{<\lambda}$. Let $\kappa < \lambda$. Let $\delta$ be a Woodin cardinal so that $\kappa < \delta < \lambda$. Let $A = p[T]$ for some $\delta^+$-weakly homogeneous tree via a $\delta^+$-weak homogeneity system $\bar\mu$. By Fact \ref{martin-solovay projection} and Fact \ref{martin-steel propogation}, $\neg A = p[\mathrm{MS}_\theta(\bar\mu)]$ and $\mathrm{MS}_\theta(\bar\mu)$ is $\kappa$-homogeneous.

Let $A \subseteq \bairespace \times \bairespace$ be in $\text{Hom}_{<\lambda}$. Let $\kappa < \lambda$. Let $\delta$ be a Woodin cardinal so that $\kappa < \delta < \lambda$. By Fact \ref{weakly homogeneous and homogeneous sets}, $\exists^\bbR A$ is $\delta^+$-weakly homogeneously Suslin via a $\delta^+$-weak homogeneity system $\bar\mu$. By Fact \ref{martin-solovay projection} and Fact \ref{martin-steel propogation}, $\mathrm{MS}_\theta(\bar\mu)$ is $\kappa$-homogeneously Suslin and $\forall^\bbR A = \neg \exists^\bbR A = p[\mathrm{MS}_\theta(\bar\mu)]$. 
\end{proof}

\Begin{fact}{coanalytic sets homogeneously suslin}
(Martin) If $\kappa$ is a measurable cardinal, then every $\coanalytic$ set is $\kappa$-homogeneously Suslin.
\end{fact}

\begin{proof}
See \cite{Determinacy-in-L(R)}, Theorem 4.15.
\end{proof}

\Begin{fact}{projective set homogeneously suslin}
(Martin-Steel) Let $\lambda$ be a limit of Woodin cardinals, then all projective sets are in $\mathrm{Hom}_{<\lambda}$. 
\end{fact}

\begin{proof}
Every Woodin cardinal has a stationary set of measurable cardinals below it. Hence every $\coanalytic$ set is $\kappa$-homogeneously Suslin for all $\kappa < \lambda$. That is, all $\coanalytic$ sets are in $\mathrm{Hom}_{<\lambda}$. Then by closure given by Fact \ref{closure of Hom}, all projective sets are in $\mathrm{Hom}_{<\lambda}$. 
\end{proof}

In fact, an even larger class of sets of reals can be homogeneously Suslin: $L(\bbR)$ is the smallest transitive class model of $\mathsf{ZF}$ containing all the reals of $V$, i.e. $(\bairespace)^V \subseteq L(\bbR)$. 

\Begin{fact}{L(R) homogeneously Suslin}
(Woodin) Suppose $\lambda$ is a limit of Woodin cardinals and there is a measurable cardinal greater than $\lambda$. Then every subsets of $\bairespace$ in $L(\bbR)$ is in $\mathrm{Hom}_{<\lambda}$. 
\end{fact}

In the previous section, sets given by projections of certains trees were essentially identified with their trees. Homogeneously Suslin sets were defined to be those sets that can be presented as projections of some trees satisfying certain properties. In the ground model, there could be many homogeneous trees representing the same homogeneously Suslin set $A$. When considering generic extensions of the ground model, there is a question of which tree should be  used to represent $A$ in the generic extension. For instance, suppose $\kappa_1 < \kappa_2$. In the ground model, suppose $A = p[T_1]$ where $T_1$ is a $\kappa_1$-homogeneous tree and $A = p[T_2]$ where $T_2$ is a $\kappa_2$-homogenous tree. Suppose $\bbP_1$ and $\bbP_2$ are two different forcing. Which tree should represent $A$ in each forcing extension? Are there circumstances in which one tree may be preferable over another? What are the relations between $p[T_1]$ and $p[T_2]$ in various forcing extensions?

Absolutely complemented trees and universally Baireness provide a way to interpreting homogeneously Suslin sets in a way which is independent of the homogeneous tree representation in some sense:

\Begin{definition}{universally baire sets}
(See \cite{Universally-Baire-Sets-of-Reals}) Let $\kappa$ be an ordinal. Let $T$ be a tree on $\omega \times X$ and let $U$ be a tree on $\omega \times Y$, for some sets $X$ and $Y$. $T$ and $U$ are $\kappa$-absolute complements if and only if for all forcings $\bbP \in V_\kappa$ and all $G \subseteq \bbP$ which are $\bbP$-generic over $V$, $V[G] \models p[T] = \bairespace \setminus p[U]$.

A tree $T$ on $\omega \times X$ is $\kappa$-absolutely complemented if and only if there exists some tree $U$ on $\omega \times Y$ (for some set $Y$) so that $T$ and $U$ are $\kappa$-absolute complements.  

A set $A \subseteq \bairespace$ is $\kappa$-universally Baire if and only if $A = p[T]$ for some tree $T$ which is $\kappa$-absolutely complemented.
\end{definition}

\Begin{fact}{independence of tree rep}
Let $T_1$ and $T_2$ be trees on $\omega \times \gamma_1$ and $\omega \times \gamma_2$ which are $\kappa$-absolutely complemented. If $\bbP \in V_\kappa$ and $G \subseteq \bbP$ is $\bbP$-generic over $V$, then $V[G] \models p[T_1] = p[T_2]$. 
\end{fact}

\begin{proof}
Let $U_1$ and $U_2$ be trees witnessing that $T_1$ and $T_2$ are $\kappa$-absolutely complemented, respectively. Suppose without loss of generality that $V[G] \models p[T_1] \cap (\bairespace \setminus p[T_2]) \neq \emptyset$. Since $T_2$ and $S_2$ are $\kappa$-absolutely complementing, $V[G] \models p[T_1] \cap p[S_2] \neq \emptyset$. Define a tree $T_1 \otimes S_2$ by
$$(s,h,g) \in T_1 \otimes S_2 \Leftrightarrow (s,h) \in T_1 \wedge (s,g) \in S_2$$
In $V[G]$, $T_1 \otimes S_2$ is ill-founded. By the absoluteness of well-foundedness, $V \models T_1 \otimes S_2$ is ill-founded. So $V \models p[T_1] \cap p[S_2] \neq \emptyset$. This is impossible since in $V$, $p[T_1] = p[T_2]$, $p[S_1] = p[S_2]$, $p[T_1] = \bairespace \setminus p[S_1]$, and $p[T_2] = \bairespace \setminus p[S_2]$. 
\end{proof}

So if $A$ is a $\kappa$-universally Baire set and if $T_1$ and $T_2$ are two $\kappa$-absolutely complemented trees so that $V \models A = p[T_1] = p[T_2]$, then either tree can be used to represent $A$ in forcing extensions by forcings in $V_\kappa$.  As a matter of convention, if $A$ is $\kappa$-universally Baire and $\bbP \in V_\kappa$, the set $A$ will always refer to $p[T]$ for some and any $\kappa$-absolutely complemented tree $T \in V$ so that $V \models p[T] = A$. 

\Begin{fact}{homogenous sets are universally baire}
Let $\kappa$ be a cardinal. $\kappa$-weakly homogenously set are $\kappa$-universally Baire.
\end{fact}

\begin{proof}
(See \cite{The-Derived-Model-Theorem}, Corollary 1.21) Let $A = p[T]$ where $T$ is a $\kappa$-weakly homogeneously Suslin set via $\kappa$-weak homogeneity system $\bar\mu$. Fact \ref{martin-solovay generically correct represention of complement} implies that for an appropriate $\theta$, $\mathrm{MS}_\theta(\bar\mu)$ witnesses that $T$ is $\kappa$-absolutely complemented.
\end{proof}

In particular, $\kappa$-homogeneously Suslin sets can be interpreted unambiguously in $\bbP$-extensions whenever $\bbP \in V_\kappa$. 

Let $\lambda$ be a limit of Woodin cardinals. Let $\dot A$ be a new unary relation symbol. Let $A \subseteq (\bairespace)^n$ be such that $A \in \text{Hom}_{<\lambda}$. $(H_{\aleph_1}, \in, A)$ be the $\{\dot\in,\dot A\}$-structure with domain $H_{\aleph_1}$ (the hereditarily countable sets) and $\dot A$ interpreted as $A$. Now let $\bbP \in V_\lambda$ be some forcing and $G \subseteq \bbP$ be a $\bbP$-generic filter over $V$. $\bbP \in V_\kappa$ for some $\kappa < \lambda$. The structure $(H_{\aleph_1}, \in, A^{V[G]})$ is understood in the following way: It is a structure with domain $H_{\aleph_1}^{V[G]}$ (the hereditarily countable subsets of $V[G]$) and $A^{V[G]}$ is $p[T]^{V[G]}$ for any $\gamma$-homogeneous tree $T$ so that $V \models A = p[T]$ and $\gamma \geq \kappa$. By the above discussion, this is independent of which tree $T$ is chosen. Actually, in the proof of the fact below, depending on the quantifier complexity of a particular formula $\varphi$ involving $\dot A$, $A$ will be considered as $p[T]$ for a sufficiently homogeneous tree $T$ so that after the appropriate number of applications of the Martin-Solovay tree construction, the resulting tree representation of $\varphi$ will be at least $\kappa$-homogeneous.

Using ideas very similar to the proof of Fact \ref{closure of Hom} (also see the proof of Fact \ref{AD absoluteness for formula real quantification over projection of tree} for Cohen forcing), one has the following absoluteness result:

\Begin{fact}{homogeneously suslin absoluteness}
(Woodin) Let $\lambda$ be a limit of Woodin cardinals. Let $A \in \text{Hom}_{<\lambda}$. Let $\bbP \in V_\lambda$ and $G \subseteq \bbP$ be $\bbP$-generic over $V$. Then $(H_{\aleph_1}^V, \in, A)$ and $(H_{\aleph_1}^{V[G]}, \in, A^{V[G]})$ are elementarily equivalent. 
\end{fact}

\begin{proof}
See \cite{The-Derived-Model-Theorem}, Theorem 2.6. 
\end{proof}

In this setting, $V$ and $V[G]$ satisfy the same formulas involving $\dot A$ and quantifications over the reals with the above intended interpretation. In particular, $V$ and $V[G]$ satisfy the same projective formulas. 

Now, the above discussion will be applied to indicate when assumption $\mathsf{A_\Sigma}$, $\mathsf{B_\Sigma}$, $\mathsf{C_\Sigma}$, $\mathsf{A_\Pi}$, $\mathsf{B_\Pi}$, and $\mathsf{C_\Pi}$.

Let $\lambda$ be a limit of Woodin cardinals. By the above discussion about universal Baireness, one may speak about an equivalence relation $E \in \mathrm{Hom}_{<\lambda}$ without explicit reference to a fix tree defining $E$. By Fact \ref{closure of Hom}, if $E \in \mathrm{Hom}_{<\lambda}$, then $\bairespace \setminus E \in \mathrm{Hom}_{<\lambda}$. Given an $\kappa$-weakly homogeneous tree representation of $E$ for sufficiently large $\gamma$, the associated Martin-Solovay tree will be a sufficiently homogeneous tree representation of $\bairespace \setminus E$ by Fact \ref{martin-steel propogation}. Hence in this setting, $E_S = E^{T}$, where $T$ is the appropriate Martin-Solovay tree using the homogeneity system on $S$. (So if $E_S$ has all $\coanalytic$ classes, then the results of Section \ref{The Game} should be applied to $E^T$ using assumption $\mathsf{A_\Pi}$, $\mathsf{B_\Pi}$, $\mathsf{C_\Pi}$ for $T$ and $I$.) Fix a $\sigma$-ideal $I$ on $\bairespace$ so that $\bbP_I$ is proper.

The formula $D_\Sigma$ and $D_\Pi$ both involve complements and real quantification over the homogeneously Suslin set $E$. By Fact \ref{closure of Hom}, $D_\Sigma, D_\Pi \in \mathrm{Hom}_{<\lambda}$. Starting with an appropriate weakly homogeneous tree representation of $E$, the process described in the proof of Fact \ref{closure of Hom} produces a tree $U$ representing $D_\Sigma$ or $D_\Pi$ which is generically correct for $\bbP_I$, in the sense that $1_{\bbP_I} \forces_{\bbP_I} p[\check U] = \{(x,T) : D_\Sigma(x,T)\}$. So assumption $\mathsf{C_\Sigma}$ holds for $E$ and $I$. (A similar argument holds for $\mathsf{C_\Pi}$.)

$E$ having all $\analytic$ classes can be expressed as a formula using some real quantifiers over the equivalence relation $E \in \mathrm{Hom}_{<\lambda}$. Fact \ref{homogeneously suslin absoluteness} implies that these statements are absolute to the $\bbP_I$-extension. The tree $S$ remains homogeneous in the $\bbP_I$-extension by the remark mentioned before Fact \ref{martin-solovay generically correct represention of complement}. This shows that $\mathsf{A_\Sigma}$ and $\mathsf{B_\Sigma}$ holds for $E$ and $I$. 

Finally, using the above discussion and results of the previous section, the following can be obtained:

\Begin{theorem}{homogeneously suslin canonicalization}
Let $\lambda$ be a limit of Woodin cardinals. Let $I$ be a $\sigma$-ideal on $\bairespace$ so that $\bbP_I$ is proper. Let $E \in \mathrm{Hom}_{<\lambda}$ be an equivalence relation on $\bairespace$. 

If $E$ has all $\analytic$ ($\coanalytic$ or $\borel$) classes, then for every $I^+$ $\borel$ set $B$, there is an $I^+$ $\borel$ $C \subseteq B$ so that $E \upharpoonright C$ is $\analytic$ ($\coanalytic$ or $\borel$, respectively).
\end{theorem}

\Begin{theorem}{projective canonicalization}
Suppose there are infinitely many Woodin cardinals. Let $I$ be a $\sigma$-ideal on $\bairespace$ so that $\bbP_I$ is proper. Let $E$ be a projective equivalence relation on $\bairespace$. 

If $E$ has all $\analytic$ ($\coanalytic$, $\borel$) classes, then for every $I^+$ $\borel$ set $B$, there is an $I^+$ $\borel$ $C \subseteq B$ so that $E \upharpoonright C$ is $\analytic$ ($\coanalytic$, $\borel$, respectively).
\end{theorem}

\Begin{theorem}{L(R) canonicalization}
Suppose there is a measurable cardinal with infinitely many Woodin cardinals below it. Let $I$ be a $\sigma$-ideal on $\bairespace$ so that $\bbP_I$ is proper. Let $E \in L(\bbR)$ be an equivalence relation on $\bairespace$. 

If $E$ has all $\analytic$ ($\coanalytic$, $\borel$) classes, then for every $I^+$ $\borel$ set $B$, there is an $I^+$ $\borel$ $C \subseteq B$ so that $E \upharpoonright C$ is $\analytic$ ($\coanalytic$, $\borel$, respectively).
\end{theorem}

With the appropriate assumptions, even more sets of reals are homogeneously Suslin and these canonicalization results would hold for equivalence relations in those classes. For example, Chang's model $L({}^\omega \text{ON}) = \bigcup_{\alpha \in \text{ON}} L({}^{\omega} \alpha)$ is the smallest inner model of $\mathsf{ZF}$ containing all the countable sequences of ordinals of $V$. Woodin has shown that with a proper class of Woodin cardinals, every set of reals in $L({}^\omega\text{ON})$ is $\infty$-homogeneously Suslin. Hence under this assumption, the above result would hold for equivalence relations in $L({}^\omega\text{ON})$ with all $\analytic$, $\coanalytic$, or $\borel$ classes.

All the above theorems also hold for graphs $G$ so that for all $x$, $G_x$ is $\analytic$ ($\coanalytic$ or $\borel$). (See the end of Section \ref{The Game}.)

\section{Canonicalization for All Equivalence Relations}\label{Canonicalization for All Equivalence Relations}
This section will consider Question \ref{every equivalence relation canonicalization}: Is it consistent that for every equivalence relation $E$ with all $\borel$ classes and every $\sigma$-ideal $I$ such that $\bbP_I$ is proper, there is an $I^+$ $\borel$ subset $C$ such that $E \upharpoonright C$ is a $\borel$ equivalence relation.

As with other regularity properties, this question has a negative answer if the axiom of choice holds. First, a definition and a property of all $\coanalytic$ equivalence relations:

\Begin{definition}{thin equivalence relation}
An equivalence relation $E$ on $\bairespace$ is thin if and only if there does not exists a perfect set $P \subseteq \bairespace$ such that $\neg(x \ E \ y)$ for all $x,y \in P$ with $x \neq y$. 
\end{definition}

There are $\analytic$ thin equivalence relation with uncountably many classes. In fact, there are $\analytic$ thin equivalence relation with all $\borel$ classes and uncountably many classes: for example, the countable admissible ordinal equivalence relation, $F_{\omega_1}$, and any counterexamples to Vaught's conjecture (if they exist). The Silver's dichotomy imply that there are no $\coanalytic$ thin equivalence relations:

\Begin{fact}{silver's dichotomy}
(Silver) If $E$ is a $\coanalytic$ equivalence relation on $\bairespace$, then either $E$ has countably many classes or there exists a perfect set of pairwise $E$-inequivalent elements. 
\end{fact}

\begin{proof}
See \cite{Counting-the-Number-of-Equivalence-Classes}.
\end{proof}

\Begin{proposition}{failure of canonicalization with axiom of choice}
$(\mathsf{ZF})$ If there is a well-ordering of $\bairespace$, then there is a thin equivalence relation $E^*$ on $\bairespace$ with equivalence classes of size at most two. 

For any $\sigma$-ideal $I$ on $\bairespace$ and any $I^+$ $\borel$ set $C$, $E^* \upharpoonright C$ is not $\borel$.
\end{proposition}

\begin{proof}
First a remark: Proposition \ref{negative answer delta12 in L} is proved in a similar way by showing that in $L$, there is a thin $\Delta_2^1$ equivalence relation with all countable classes. 

Now the proof of the proposition: Using the well-ordering of $\bairespace$, let $\Phi : 2^{\aleph_0} \rightarrow \bairespace$ be bijection and let $\Psi : 2^{\aleph_0} \rightarrow \bairespace$ be an enumeration of all the perfects trees on $\omega$. 

The equivalence $E^*$ is defined by stages through transfinite recursion as follows: 

Let $A_0 = \emptyset$. $E^*_0 = \emptyset$.

Stage $\xi + 1$: Suppose $A_\xi$ and $E^*_\xi$ have been defined with $|A_\xi| < 2^{\aleph_0}$. Find some reals $r_\xi$ and $s_\xi$ so that $r_\xi, s_\xi \notin A_\xi$, $r_\xi \neq s_\xi$, and $r_\xi, s_\xi \in [\Psi(\xi)]$. 

If $\Phi(\xi) \in A_\xi \cup \{r_\xi, s_\xi\}$, then define $A_{\xi + 1} = A_\xi \cup \{r_\xi, s_\xi\}$ and 
$$E^*_{\xi + 1} = E^*_\xi \cup \{(r_\xi,r_\xi), (s_\xi,s_\xi), (r_\xi, s_\xi), (s_\xi,r_\xi)\}$$
If $\Phi(\xi) \notin A_\xi \cup \{r_\xi, s_\xi\}$, then define $A_{\xi + 1} = A_\xi \cup \{r_\xi, s_\xi, \Phi(\xi)\}$ and 
$$E^*_{\xi + 1} = E^*_\xi \cup \{(\Phi(\xi), \Phi(\xi)), (r_\xi,r_\xi), (s_\xi,s_\xi), (r_\xi, s_\xi), (s_\xi,r_\xi)\}$$

At limit stage $\xi$: Let $A_\xi = \bigcup_{\eta < \xi} A_\eta$ and $E_\xi^* = \bigcup_{\eta < \xi} E^*_\eta$. 

Note that $A_{2^{\aleph_0}} = \bairespace$. Let $E^* = E^*_{2^{\aleph_0}}$. $E^*$ is an equivalence relation on $\bairespace$. $E^*$ has classes of size at most two. $E^*$ is thin: Suppose $T$ is a perfect tree on $\omega$. Then $T = \Psi(\xi)$ for some $\xi < 2^{\aleph_0}$. Then $r_\xi \ E^* \ s_\xi$ and $r_\xi,s_\xi \in [\Psi(\xi)] = [T]$. 

Now let $I$ be a $\sigma$-ideal on $\bairespace$. Suppose there was some $I^+$ $\borel$ set $C$ so that $E^* \upharpoonright C$ is $\borel$. Since $C$ is $I^+$ and $I$ is a $\sigma$-ideal, $C$ must be uncountable. Since $E^*$ has classes of size at most two, $E \upharpoonright C$  can not have only countably many classes. Since $\borel$ equivalence relations are $\coanalytic$, the Silver's dichotomy (Fact \ref{silver's dichotomy}) implies that there is a perfect set $P \subseteq C$ of $E^*$-inequivalent elements. There is a perfect tree $T$ so that $[T] = P$. Let $\xi < 2^{\aleph_0}$ be so that $\Psi(\xi) = T$. Then $r_\xi, s_\xi \in [T] = P \subseteq C$ and $r_\xi \ E^* \ s_\xi$. Contradiction.
\end{proof}

Hence to get a positive answer to Question \ref{every equivalence relation canonicalization}, there can not exist a well-ordering of the reals, so the full axiom of choice must fail.

First, the immediate concern in the choiceless setting is the definition of properness: Since set may not have a cardinality, it is preferable to use $V_\Xi$ rather than $H_\Xi$. Recall in $\mathsf{ZFC}$, for any $\sigma$-ideal $I$ on $\bairespace$, $\bbP_I$ was proper if and only if for all sufficiently large cardinals $\Xi$, any $B \in \bbP_I$, and all countable elementary $M \prec V_\Xi$ with $\bbP_I, B \in M$, $\{x \in B : x \text{ is } \bbP_I\text{-generic over } M\}$ is $I^+$ $\borel$. Without the axiom of choice, the downward Lowenheim-Skolem theorem may fail for structures in countable languages and so there may be no countable elementary substructure. Moreover, in the previous section, it was also important to be able to choose countable elementary substructures containing certain homogeneously Suslin trees. 

However, only dependence choice ($\mathsf{DC}$) is needed to prove the following form of the downward Lowenheim-Skolem theorem: Let $\SCRL$ be a countable language. Let $M$ be an $\SCRL$-structure. Let $A \subseteq M$ be countable. Then there exists an $\SCRL$-elementary substructure $N$ of $M$ so that $A \subseteq N$. 

Hence with $\mathsf{DC}$, the definition of properness and the ability to construct elementary substructure of $V_\Xi$ with certain desired objects inside are still available. 

Without the axiom of choice, determinacy for various games are useful for settling many questions in descriptive set theory: The axiom of determinacy ($\mathsf{AD}$) asserts that all games of the form in Definition \ref{games} where the moves are elements of $\omega$ are determined. The axiom of determinacy for the reals ($\mathsf{AD_\bbR}$) asserts that all games of the form in Definition \ref{games} where the moves are elements of $\bairespace$ are determined. 

In terms of large cardinals, the consistency of $\mathsf{AD}$ follows from the consisteny of infinitely many Woodin cardinals. The consistency of $\mathsf{AD_\bbR}$ follows from the consistency of the existence of a cardinal $\lambda$ which is both a limit of Woodin cardinals and $<\lambda$-strong cardinals. 

$\Theta$ denotes the supremum of the ordinals which are surjective image of $\bbR$. As described above, $\mathsf{DC}$ would be useful for carrying out arguments from the earlier sections. A result of Solovay shows that $\mathsf{ZF} + \mathsf{AD_\bbR} + V = L(\mathscr{P}(\bbR)) + \mathrm{cof}(\Theta) > \omega$ can prove $\mathsf{DC}$. It should be noted that Solovay has also shown that $\mathsf{ZF} + \mathsf{AD_\bbR} + \mathsf{DC}$ can prove the consistency of $\mathsf{ZF} + \mathsf{AD_\bbR}$; hence $\mathsf{AD_\bbR} + \mathsf{DC}$ is strictly stronger than $\mathsf{AD_\bbR}$ in terms in consistency. (See \cite{The-Independence-of-DC-from-AD} for these results concerning $\mathsf{AD_\bbR}$ and $\mathsf{DC}$.)

$\mathsf{AD_\bbR}$ is preferable over $\mathsf{AD}$ since $\mathsf{AD_\bbR}$ can prove that every subset of $\bairespace$ is homogeneously Suslin and can prove a strong form of absoluteness for proper forcings: 

\Begin{fact}{AD(R) every tree weakly homogeneously suslin}
(Martin) Under $\mathsf{ZF} + \mathsf{AD_\bbR}$, every tree on $\omega \times \lambda$, where $\lambda$ is an ordinal, is weakly homogeneously Suslin.
\end{fact}

\begin{proof}
See \cite{Weakly-Homogeneous-Trees}.
\end{proof}

\Begin{fact}{AD suslin cosuslin homogeneous suslin}
(Martin) Under $\mathsf{ZF} + \mathsf{DC} + \mathsf{AD}$, for every $A \subseteq \bairespace$, $A$ is homogeneously Suslin if and only if if $A$ and $\bairespace \setminus A$ are Suslin. Moreover, one can find a homogeneously Suslin tree $T$ on $\omega \times \kappa$, for $\kappa < \Theta$, so that $A = p[T]$. 
\end{fact}

\begin{proof}
See \cite{Tree-of-a-Moschovakis-Scale}. 
\end{proof}

\Begin{fact}{AD(R) and suslin cosuslin}
(Martin, Woodin) Under $\mathsf{ZF} + \mathsf{AD}$, $\mathsf{AD_\bbR}$ is equivalent to the statement that every subset of $\bairespace$ is Suslin.
\end{fact}

Combining the last two facts give:

\Begin{fact}{AD(R) every set is homogeneously Suslin}
Under $\mathsf{ZF} + \mathsf{AD_\bbR}$, every subset of the $\bairespace$ is homogeneously Suslin.
\end{fact}

In the previous section, an important aspect of analyzing tree representations in generic extensions was the fact that any $\kappa$-complete measure $\mu$ could be naturally extended to a $\kappa$-complete measure $\mu^*$ in a forcing extension by $\bbP$, whenever $|\bbP| < \kappa$. 

Let $I$ be a $\sigma$-ideal on $\bairespace$. $\bbP_I$ is in bijection with $\bairespace$ and hence is not well-ordered under $\mathsf{AD}$. Also note that the measures produced using $\mathsf{AD}$ to witness homogeneity and weak homogeneity are $\aleph_1$-complete. For the general $\sigma$-ideal, it is not clear how to modify the arguments of the previous section in the context of $\mathsf{AD}_\bbR$.

However, there is one important $\sigma$-ideal for which the previous arguments will work with minor modifications: For the meager ideal, $\bbP_{I_\mathrm{meager}}$ is forcing equivalent to Cohen forcing, denoted $\bbC$, which is a countable forcing. 

Let $T$ be an $\aleph_1$-weakly homogeneous tree on $\omega \times \gamma$ witnessed by the weak homogeneity system $\bar{\mu}$. Fact \ref{martin-solovay projection}, which is provable in $\mathsf{ZF} + \mathsf{DC}$, implies that $V \models p[T] = \bairespace \setminus \mathrm{MS}_{\gamma^+}(\bar{u})$. 

Since $|\bbC| = \aleph_0 < \aleph_1$, any $\aleph_1$-complete measure can be extended to an $\aleph_1$-complete measure in the $\bbC$-forcing extension. Likewise, every $\aleph_1$-weak homogeneity system $\bar{\mu}$ can be extended to an $\aleph_1$-weak homogeneity system. Fact \ref{martin-solovay generically correct represention of complement} and the discussion before it holds when $\kappa = \aleph_1$ and $\bbP = \bbC$: 

\Begin{fact}{martin-solovay projection cohen AD(R)}
Assume $\mathsf{ZF} + \mathsf{AD}_\bbR$. Let $T$ be an $\aleph_1$-weakly homogeneous tree on $\omega \times \gamma$ witnessed by the weak homogeneity system $\bar{\mu}$. If $G \subseteq \bbC$ is $\bbC$-generic over $V$, then $V[G] \models \mathrm{MS}_{\gamma^+}(\bar{\mu}^V) = \mathrm{MS}_{\gamma^+}(\bar{\mu}^*)$ and $V[G] \models p[\mathrm{MS}_{\gamma^+}(\bar{\mu})^V] = \bairespace \setminus p[T]$. 
\end{fact}

To show that $T$ and $\mathrm{MS}_{\gamma^+}(\bar{\mu})$ continues to complement each other in the $\bbP_I$ generic extension for an arbitrary $\sigma$-ideal $I$ so that $\bbP_I$ is proper, a strong absoluteness result for proper forcing due to Neeman and Norwood will be used. This result is similar to \cite{Proper-Forcing-and-Absoluteness}. 

\Begin{fact}{AD(R) embedding theorem}
(Neeman and Norwood, to appear) Under $\mathsf{ZF} + \mathsf{AD}_\bbR + V = L(\mathscr{P}(\bbR))$, for every proper forcing $\bbP$, and $G \subseteq \bbP$ which is $\bbP$-generic over $V$, there is an elementary embedding $j : L(\mathscr{P}(\bbR)) \rightarrow L(\mathscr{P}(\bbR)^{V[G]})$ so that $j$ does not move ordinals or reals. 
\end{fact}

\Begin{fact}{tree projection AD(R) + V = L(P(R)) proper forcing}
Assume $\mathsf{ZF} + \mathsf{AD}_\bbR + V = L(\mathscr{P}(\bbR))$. Suppose $T$ and $S$ are trees on $\omega \times \gamma$ and $\omega \times \delta$ so that $p[T] = \bairespace \setminus p[S]$. Then $V[G] \models p[T] = \bairespace \setminus p[S]$.

In particular, if $T$ is a weakly homogeneous tree on $\omega \times \gamma$ witnessed by the weak homogeneity system $\bar{\mu}$, then $V[G] \models p[T] = \bairespace \setminus p[\mathrm{MS}_{\gamma^+}(\bar{\mu})]$. 
\end{fact}

\begin{proof}
Let $G \subseteq \bbP$ be $\bbP$-generic over $V$. Let $j : L(\mathscr{P}(\bbR)) \rightarrow L(\mathscr{P}(\bbR)^{V[G]})$ be an elementary embedding which does not move ordinals or reals. Note that if $T$ is a tree on $\omega \times \gamma$, then $j(T)$ is a tree on $j(\omega) \times j(\gamma) = \omega \times \gamma$ and for all $s \in {}^{<\omega}(\omega \times \gamma)$, $s \in T$ if and only if $j(s) \in j(T)$ if and only if $s \in j(T)$. Hence $T = j(T)$ and similarly $S = j(S)$. So by elementarity, $L(\mathscr{P}(\bbR)^{V[G]}) \models p[T] = \bairespace \setminus p[S]$. As $V[G]$ and $L(\mathscr{P}(\bbR)^{V[G]})$ have the same reals, $V[G] \models p[T] = \bairespace \setminus p[S]$.

For the second statement, note that under $\mathsf{ZF} + \mathsf{DC}$, Fact \ref{martin-solovay projection} implies that $L(\mathscr{P}(\bbR)) \models p[T] = \bairespace \setminus p[\mathrm{MS}_{\gamma^+}(\bar{\mu})]$. The rest follows by applying the first part.
\end{proof}

\Begin{fact}{AD real quantifier over tree representation has tree representation}
Assume $\mathsf{ZF} + \mathsf{DC} + \mathsf{AD}_\bbR + V = L(\mathscr{P}(\bbR))$ (or just $\mathsf{ZF} + \mathsf{DC} + \mathsf{AD}_\bbR$ for Cohen forcing, $\bbC$). Let $\bbP$ be a proper forcing. Suppose $T$ is a tree on ${}^k\omega \times \gamma$ for some cardinal $\gamma$ and $k \in \omega$. If $A \subseteq {}^j(\bairespace)$, for some $j \leq k$, is defined by applying complementation and $\exists^\bbR$ over $p[T]$, then there is some tree $U$ on ${}^j\omega \times \delta$ for some cardinal $\delta$ so that $A = p[U]$ and $1_\bbP \forces_\bbP A = p[\check U]$.
\end{fact}

\begin{proof}
This is proved by induction. Suppose $B$ is some set defined by real quantifiers over $p[T]$ such that there is some tree $L$ on ${}^l\omega\times \epsilon$ so that $B = p[L]$ and $1_\bbP \forces_\bbP B = p[L]$.

For the $\exists^\bbR$ case: Suppose $l = i + 1$. Define the tree $U$ on ${}^i \omega \times \epsilon$ as the induced tree defined by considering the tree $L$ on ${}^{i + 1}\omega \times \omega$ as a tree on ${}^i\omega \times (\omega \times \epsilon)$ with $\epsilon$ and $\omega \times \epsilon$ identified by some bijection. Then $\exists^\bbR B = p[U]$. 

For complementation: By Fact \ref{AD(R) every tree weakly homogeneously suslin}, $L$ is weakly homogeneously Suslin. Let $\bar{\mu}$ be some weak homogeneity system witnessing this for $L$. By Fact \ref{martin-solovay projection}, $p[T] = \bairespace \setminus p[\mathrm{MS}_{\epsilon^+}(\bar{\mu})]$. By Fact \ref{tree projection AD(R) + V = L(P(R)) proper forcing} (or Fact \ref{martin-solovay projection cohen AD(R)}), $1_\bbP \forces_\bbP p[T] = \bairespace \setminus p[\mathrm{MS}_{\epsilon^+}(\bar{\mu})]$.
\end{proof}

\Begin{fact}{AD absoluteness for formula real quantification over projection of tree}
Assume $\mathsf{ZF} + \mathsf{DC} + \mathsf{AD}_\bbR + V = L(\mathscr{P}(\bbR))$ (or just $\mathsf{ZF} + \mathsf{DC} + \mathsf{AD}_\bbR$ in the case of $\bbC$). Let $\bbP$ be a proper forcing. Let $T$ be a tree on ${}^k\omega \times \gamma$ for some cardinal $\gamma$. Let $A$ denote a predicate symbol for $p[T]$ which will always be interpreted as $p[T]$ in forcing extensions. Let $\varphi$ be a formula on $\bbR$ using predicate $A$, complementation, and $\exists^\bbR$. Then for all $r \in \bbR^V$, $V \models \varphi(r) \Leftrightarrow V[G] \models \varphi(r)$, whenever $G \subseteq \bbP$ is $\bbP$-generic over $V$.
\end{fact}

\begin{proof}
In the $V = L(\mathscr{P}(\bbR))$ case, this is essentially immediate from the absoluteness result of Fact \ref{AD(R) embedding theorem} and the fact that $L(\mathscr{P}(\bbR)^{V[G]}) \models \varphi(r) \Leftrightarrow V[G] \models \varphi(r)$.

So consider the case for $\bbC$: Let $G \subseteq \bbC$ be a $\bbC$-generic over $V$. By Fact \ref{AD real quantifier over tree representation has tree representation}, for some tree $U$, $V \models (\forall x)(\varphi(x) \Leftrightarrow x \in p[U])$ and $V[G] \models (\forall x)(\varphi(x) \Leftrightarrow x \in p[U])$. Then for any $x \in \bbR^{V}$,
$$V \models \varphi(x) \Leftrightarrow V \models x \in p[U] \Leftrightarrow V[G] \models x \in p[U] \Leftrightarrow V[G] \models \varphi(x)$$
where the second equivalence follows from the absoluteness of well-foundedness.
\end{proof}

Now assume $\mathsf{ZF} + \mathsf{DC} + \mathsf{AD}_\bbR + V = L(\mathscr{P}(\bbR))$ (or just $\mathsf{ZF} + \mathsf{DC} + \mathsf{AD}_\bbR$ when working with the meager ideal). Let $E$ be an equivalence relation on $\bairespace$ with all $\analytic$ (or $\coanalytic$ classes). Let $I$ be a $\sigma$-ideal on $\bairespace$ so that the associated forcing $\bbP_I$ is a proper forcing.

By Fact \ref{AD(R) every set is homogeneously Suslin}, $E$ is homogeneously Suslin. Let $S$ be a homogeneously Suslin tree so that $E = p[S]$. In the case of the meager ideal and under $\mathsf{ZF} + \mathsf{DC} + \mathsf{AD}_\bbR$: By the countability of $\bbC$, the argument above about weak homogeneity system would show that the homogeneity system for $S$ would lift to a homogeneity system for $S$ in the $\bbC$-extension. Thus $S$ would still be a homogeneous tree in the $\bbC$-extension. Under $\mathsf{ZF} + \mathsf{DC} + \mathsf{AD}_\bbR + V = L(\mathscr{P}(\bbR))$, for the general $\sigma$-ideal $I$ with $\bbP_I$ proper, Fact \ref{AD(R) embedding theorem} gives an elementary embedding $j : L(\mathscr{P}(\bbR)) \rightarrow L(\mathscr{P}(\bbR)^{V[G]})$. So $L(\mathscr{P}(\bbR)^{V[G]}) \models S$ is homogeneously Suslin. This is not exactly the requirement of $\mathsf{A_\Sigma}$ or $\mathsf{A_\Pi}$ so the proof of Theorem \ref{homogeneous suslin equiv with all analytic class analytic somewhere} needs to be slightly modified: To prove Claim 2, first use the same argument with the fact that $S$ is a homogeneous tree in $L(\mathscr{P}(\bbR)^{M[g]})$ to show that Player 2 has a winning strategy in $L(\mathscr{P}(\bbR)^{M[g]})$. This strategy is still a winning strategy for Player 2 in $M[g]$. This proves Claim 2 and the rest of the argument remains unchanged.

$\mathsf{A_\Sigma}$ (and similarly $\mathsf{A_\Pi}$) holds for $S$ and $I$, except for the minor point of the previous paragraph. The formula $D_\Sigma(x,T)$ from Definition \ref{Assumption B} can be expressed as a statement involving a predicate for $p[S]$, complementation, and real quantifiers. Fact \ref{AD real quantifier over tree representation has tree representation} shows that there is some tree $U$ representing $D_\Sigma$ in $V$ and in $\bbP_I$-extensions. Statement $\mathsf{C_\Sigma}$ (and similarly $\mathsf{C_\Pi}$) holds for $S$ and $I$. The statement $(\forall x)(\exists T)D_\Sigma(x,T)$ is true in $V$ since $E$ is an equivalence relation with all $\analytic$ classes. This formula is also expressed as a statement involving a predicate for $p[S]$, complementation, and real quantifiers, so Fact \ref{AD absoluteness for formula real quantification over projection of tree} implies that this statement remains true in the $\bbP_I$-extension. $\mathsf{B_\Sigma}$ (and similarly $\mathsf{B_\Pi}$) holds for $S$ and $I$. 

As Section \ref{The Game} works in $\mathsf{ZF} + \mathsf{DC}$, the arguments of that section can be carried out in the present context, with the changes mentioned above. (Recall the discussion earlier in this section about properness and elementary substructures under $\mathsf{DC}$.) 

Since Cohen forcing satisfies the $\aleph_1$-chain condition, one can obtain more than just canonicalization on a nonmeager set but in fact on a comeager set.

Finally the following results are obtained. Again, the analogous result for graphs also hold:

\Begin{theorem}{AD(R) canonicalization for all equiv meager ideal}
Assume $\mathsf{ZF} + \mathsf{DC} + \mathsf{AD}_\bbR $. Let $E$ be an equivalence relation on $\bairespace$. If $E$ has all $\analytic$ ($\coanalytic$ or $\borel$) classes, then for every nonmeager $\borel$ set $B$, there is a $\borel$ set $C \subseteq B$ which is comeager in $B$ so that $E \upharpoonright C$ is $\analytic$ ($\coanalytic$ or $\borel$, respectively).
\end{theorem}

\Begin{theorem}{AD(R) V = L(P(R)) canonicalization for all equiv all proper ideal}
Assume $\mathsf{ZF} + \mathsf{DC} + \mathsf{AD}_\bbR + V = L(\mathscr{P}(\bbR))$. Let $I$ be a $\sigma$-ideal on $\bairespace$ so that $\bbP_I$ is proper. Let $E$ be an equivalence relation on $\bairespace$. If $E$ has all $\analytic$ ($\coanalytic$ or $\borel$) classes, then for every $I^+$ $\borel$ set $B$, there is an $I^+$ $\borel$ set $C \subseteq B$ so that $E \upharpoonright C$ is $\analytic$ ($\coanalytic$ or $\borel$, respectively).
\end{theorem}

\bibliographystyle{amsplain}
\bibliography{references}

\end{document}